\documentclass[11pt]{article}
\usepackage{xcolor}
\usepackage{amsmath}
\usepackage{amsthm}
 \usepackage{amssymb}
 \usepackage{leftidx}
\usepackage{graphicx}
\usepackage{amsfonts}
\usepackage{graphicx, color}
 \usepackage{cite}
 \usepackage{bm}
 \usepackage{leftidx}
 \numberwithin{equation}{section}

\newtheorem{theorem}{Theorem}[section]
\newtheorem{lemma}{Lemma}[section]

\newtheorem{remark}{Remark}
\renewenvironment{proof}[1][Proof]{\noindent\textit{#1. } }{\hfill$\square$}
\newcommand{\fd}[1]{\mathcal{D}^{#1}_t}
\newcommand{\dfd}[1]{\mathcal{D}^{#1}_\tau}

\newcommand{\diff}{\triangledown_\tau}
\newcommand{\Ass}[1]{\textbf{\upshape A#1}}

\newcommand{\iprod}[1]{\left\langle#1\right\rangle}
\newcommand{\taumax}{\tau}
\newcommand{\rhomax}{\rho}

\newcommand{\defeq}{:=}

\newcommand{\zd}{\,\mathrm{d}}
\newcommand{\abs}[1]{\left|#1\right|}
\newcommand{\absb}[1]{\big|#1\big|}

\newcommand{\bra}[1]{\left(#1\right)}
\newcommand{\brab}[1]{\big(#1\big)}
\newcommand{\braB}[1]{\Big(#1\Big)}

\newcommand{\mynormb}[1]{\big\|#1\big\|}

\newcommand{\cy}[1]{{\color{cyan}#1}}

\newcommand{\Caputo}{{{}^{\text{C}}_0\fd{\alpha}}}
\newcommand{\DCaputo}{{\mathcal{D}^\alpha_\tau}}

\newcommand{\Gloc}{G_{\mathrm{loc}}}
\newcommand{\Ghis}{G_{\mathrm{his}}}

\newcommand{\Mss}[1]{\textbf{\upshape M#1}}

\newcommand{\TheTitle}{$\alpha$-robust error estimates of general non-uniform time-step numerical schemes for reaction-subdiffusion problems}
\title{{\TheTitle}\thanks{This work was supported in part by NSFC grants 12171376, 12131005, 2020-JCJQ-ZD-029, and the Fundamental Research Funds for the Central Universities 2042021kf0050.}}
\author{
Jiwei Zhang\thanks{
School of Mathematics and Statistics, and Hubei Key Laboratory of Computational Science, Wuhan University, Wuhan 430072, China (jiweizhang@whu.edu.cn).}
\and Zhimin Zhang\thanks{
Department of Mathematics, Wayne State University, Detroit, MI 48202, USA (ag7761@wayne.edu).}
\and Chengchao Zhao\thanks{
Beijing Computational Science Research Center, {Beijing 100193,  China} (cheng\_chaozhao@csrc.ac.cn). }
}
\date{}

\begin{document}
\maketitle

\begin{abstract}
Numerous error estimates have been carried out on various numerical schemes for subdiffusion equations. Unfortunately most error bounds suffer from a factor $1/(1-\alpha)$ or $\Gamma(1-\alpha)$, which blows up as the fractional order $\alpha\to 1^-$, a phenomenon not consistent with regularity of the continuous problem and numerical simulations in practice. Although efforts have been made to avoid the factor  blow-up phenomenon, a robust analysis of error estimates still remains incomplete for numerical schemes with general nonuniform time steps. In this paper, we will consider the $\alpha$-robust error analysis of convolution-type schemes for subdiffusion equations with general nonuniform time-steps, and provide explicit factors in error bounds with  dependence information on $\alpha$ and temporal mesh sizes.  As illustration, we apply our abstract framework to two widely used schemes, i.e., the L1 scheme and Alikhanov's  scheme.  Our rigorous proofs reveal that the stability and convergence of a class of convolution-type schemes is $\alpha$-robust, i.e., the factor will not blowup while  $\alpha\to 1^-$ with general nonuniform time steps even when rather general initial regularity condition is considered.
\end{abstract}
{{\bf Keywords:} 
subdiffusion equations,  nonuniform time mesh, discrete Gr\"{o}nwall inequality, $\alpha$-robust error estimate
}
%\begin{AMS}
%65M06, 35B65
%\end{AMS}

\section{Introduction}
We revisit sharp error estimates of general nonuniform time-step numerical schemes for solving the following reaction-subdiffusion problem on a bounded spatial domain $\Omega$:
\begin{equation}\label{EQ_equation}
\begin{aligned}
\fd{\alpha} u- \Delta u&= \kappa u +f,&&\text{for}\; x\in\Omega, \;0<t\le T,\\
u&=u_0(x), &&\text{for}\; x\in\Omega, \;t=0,\\
u&=0, &&\text{for}\; x\in\partial\Omega, \;0<t<T,
\end{aligned}
\end{equation}
where $\kappa \in \mathbb{R}$ is the reaction coefficient and $\fd{\alpha}$ denotes the Caputo derivative by
\begin{align}
\fd{\alpha} u :=\int_0^t\omega_{1-\alpha}(t-s) u_s(s)\zd s,\quad 0<\alpha<1\label{EQ_Caputo}
\end{align}
with $\omega_\beta(t):={t^{\beta-1}}/{\Gamma(\beta)}$ and $\Gamma(\cdot)$ representing the Gamma function.

%The model \eqref{EQ_equation} is used to describe the dynamics of  anomalous diffusion, and  has various practical applications in engineering, physics and biology \cite{MK00,MK04,MJCB14}. Hence,  its effective numerical schemes and  rigorous numerical analysis are of great practical importance.

Numerous schemes have been developed over last two decades including L1 scheme \cite{SW06,LX07}, L2-type scheme \cite{LX16,CXW14,K14}, Alikhanov's scheme \cite{Alikhanov2015,LMZ21}, and convolution quadrature (CQ) method \cite{L88,JLZ18,JLZ19}. Also, there has  been an explosive growth in the  numerical analysis of these schemes.  Among the analysis of these schemes,  the error bounds  generally contain a factor $1/(1-\alpha)$, which will blow up as $\alpha\to 1^-$.  As mentioned in \cite{JLZ19} :``This phenomenon does not fully agree with the results for the continuous model ...and it is of interest to further refine the estimates to fill in the gap.''

Recently,  a refined error estimate  has been presented in \cite{CS21,HS21,HSC21,WS21} to  avoid the blow-up phenomenon,  which is the so-called $\alpha$-robust error estimate. One of the main defects in the existing analysis is  the mesh restriction for the consistence error arose from the discrete scheme of the  Caputo derivative, i.e., the analysis depends on the precise form of graded mesh $t_n = (k/N)^\gamma$ with $\gamma\geq 1$.
 On the other hand, adaptive time steps and other nonuniform temporal meshes have been tested to be practically valuable for resolving complex behaviors (such as physical oscillations, multi-scale dynamics) of linear/nonlinear
fractional differential equations. But, the $\alpha$-robust analysis of various schemes on general  nonuniform meshes still remains incomplete.

The aim of this paper is to present an $\alpha$-robust error estimate of convolution-types schemes proposed in \cite{LMZ19} for \eqref{EQ_equation} with  nonuniform time steps. Compared with traditional local methods for the first-order derivative, the numerical analysis for
nonlocal time-stepping schemes on nonuniform time meshes is challenging
due to the convolution integral (nonlocal) form of the fractional derivative. Recently, a theoretical framework \cite{LLZ18,LMZ21,LMZ19} is developed for the stability and convergence analysis of a class of widely used schemes with nonuniform meshes. The framework provides three useful tools: (i) a novel concept of discrete complementary convolution (DCC) kernels, (ii) a discrete fractional Gr\"{o}nwall inequality (DFGI) to address the stability of numerical schemes, and (iii)  a concept of \emph{error convolution structure} ({\bf ECS}) to bound the local consistency errors. Nonetheless, their error bounds contain the factor $1/(1-\alpha)$ and thus suffer from the factor blow-up phenomenon.

In this paper, we offer a new viewpoint based on aforementioned framework by analyzing how and where the error bound depends on the fractional order $\alpha$. Actually, the undesired  factor $1/(1-\alpha)$ is produced in the estimate of global consistency error for the discrete scheme of the Caputo derivative (see Remark \ref{Remark_L1} for the detail discussion). A natural observation is that the original  global consistency error is a form of summation and does not contain any blow-up factor like $1/(1-\alpha)$. Inspired by the observation, we reproduce a summation-type error bound for the global consistency error by a novel technique based on the DCC kernels and {\bf ECS}. The resulting error bound does not contain any blow-up factor and thus is $\alpha$-robust.

To illustrate the idea,
we then apply the $\alpha$-robust error estimate to widely used L1 and  Alikhanov's schemes on three types of temporal meshes. That is, uniform mesh,  the well-known  graded mesh $t_n = (k/N)^\gamma$ with $\gamma\geq 1$ and the general nonuniform mesh that the following condition hold \cite{MM15,B85}.
\begin{enumerate}
\item[\Mss{1}.] There is a constant $C_1,C_{\gamma}$ such that $\tau_1\geq C_1 \tau^\gamma$ and
$\tau_k\le C_{\gamma}\tau\min\{1,t_k^{1-1/\gamma}\}$
for~$1\le k\le N$,  $t_{k}\leq C_{\gamma}t_{k-1}$ and $\tau_k \leq  C_{\gamma} \tau_{k-1}$ for~$2\le k\le N$.
\end{enumerate}
Here the parameter $\gamma\geq 1$ controls the extent to which the time levels are concentrated near $t=0$. As $\gamma$ increases, the initial step sizes become smaller than the later ones.
The mesh \Mss{1} is also employed to reproduce the optimal convergence order in parabolic equations with nonsmooth initial values \cite{LM21,LM22,LMW21}.

Our analysis reveals that the factor $\chi_{n,\gamma}^L$ defined by (\ref{EQ_cl}) for L1 scheme and $\chi_{n,\gamma}^A$ (\ref{EQ_ca}) for Alikhanov's schemes in the error bounds will not blow up as $\alpha \to 1^{-}$. Specifically, the formulas of the factors, say $\chi_{n,\gamma}$ by removing the index $A$ or $L$, are listed  in Table \ref{Table1} for $\sigma\in (0,1)$ and in Table \ref{Table2}  for $\sigma\in (1,2)$. Here $\sigma$ is a regularity parameter appeared in the general regularity condition \eqref{EQ_2153}.  One can see from Tables \ref{Table1} and \ref{Table2} that the factors $\chi_{n,\gamma}$ are robust for all $\sigma$ and $\alpha$, i.e.,  it will not blow up as $\alpha\to 1^-$. For $\sigma\in (1,2)$, Table \ref{Table2} shows that the factor  will be $\ln n$ with  $n$ representing the index of time level or $\ln n^*$ with $n^*:= (\frac{t_n}{\tau_1})^\frac1{\gamma}$ when  $\sigma \gamma\geq2-\alpha $ for the L1 scheme  or $\sigma \gamma\geq3-\alpha$ for Alikhanov's scheme with the grading parameter $\gamma\geq 1$. In this situation, the quasi-optimal convergence  is  achieved. For $\sigma\in (0,1)$, Table \ref{Table1} shows that the factor  $\chi_{n,\gamma}$ is $\mathcal{O}(1)$ even though $\alpha\to 1^-$ excluding the special case $\sigma \gamma=2-\alpha $ for the L1 scheme  or $\sigma \gamma=3-\alpha$ for Alikhanov's scheme. In this situation, the optimal convergence  is  achieved.
For the case of $\alpha\to 1^-$, the compatible factors $\chi_{n,\gamma}$  can be bounded by $\ln n$, which will produce the uniformly quasi-optimal $\alpha$-robust error estimates for a class of widely used schemes since the error bound may appear a factor of $\ln n$.

\begin{table}[htb]
\begin{center}
\def\temptablewidth{1.0\textwidth}
{\rule{\temptablewidth}{0.7pt}}
\begin{tabular*}{\temptablewidth}{@{\extracolsep{\fill}}llll}
 %&\multicolumn{3}{c}{L1 scheme }&\multicolumn{3}{c}{Alikhanov's scheme}  \\
% \hline
L1 scheme &$\sigma \gamma<2-\alpha$& $\sigma \gamma=2-\alpha$ &$\sigma \gamma>2-\alpha$ \\
\hline
\quad Uniform mesh  &$\frac{1-n^{\sigma+\alpha-2}}{2-\alpha -\sigma}$   & $\ln n $    &$\frac{1-n^{2-\alpha-\sigma}}{\alpha-2+\sigma  }  $   \\
\quad Graded mesh  &$\frac{1-n^{\sigma\gamma+\alpha-2}}{2-\alpha -\sigma\gamma}$   & $\ln n $    &$\frac{1-n^{2-\alpha-\sigma\gamma}}{\alpha-2+\sigma\gamma}  $   \\

\quad \Mss{1} &$\frac{1-(n^*)^{\sigma\gamma+\alpha-2}}{2-\alpha -\sigma\gamma}$     &$\ln n^*$       &$\frac{1-(n^*)^{2-\alpha-\sigma\gamma}}{\alpha-2+\sigma\gamma}$          \\
% general nonuniform mesh  ($\gamma\geq 1$)   &$\ln(t_n/\tau_1) $      &$\ln(t_n/\tau_1) $      &$\ln(t_n/\tau_1) $       \\
\hline
\hline
Alikhanov's scheme &$\sigma \gamma<3-\alpha$& $\sigma \gamma=3-\alpha$ &$\sigma \gamma>3-\alpha$ \\
\hline
\quad Uniform mesh &$\frac{1-n^{\sigma+\alpha-3}}{3-\alpha -\sigma}$   &  $\ln n $     &$\frac{1-n^{3-\alpha-\sigma}}{\alpha-3+\sigma}  $    \\
\quad Graded mesh   &$\frac{1-n^{\sigma\gamma+\alpha-3}}{3-\alpha -\sigma\gamma}$   &  $\ln n $     &$\frac{1-n^{3-\alpha-\sigma\gamma}}{\alpha-3+\sigma\gamma}  $    \\
\quad \Mss{1}   &$\frac{1-(n^*)^{\sigma\gamma+\alpha-3}}{3-\alpha -\sigma\gamma}$     &$\ln n^*$      &$\frac{1-(n^*)^{3-\alpha-\sigma\gamma}}{\alpha-3+\sigma\gamma}$          \\\end{tabular*}
{\rule{\temptablewidth}{0.7pt}}
\end{center}
\caption {$\sigma\in (0,1)$:The formula of factor $\chi_{n,\gamma}^L$ for L1  and $\chi_{n,\gamma}^A
$ Alikhanov's schemes with respect to parameters of fractional order $\alpha$, regularity parameter $\sigma$, the grading parameter $\gamma$ and time step $n$.   \label{Table1}}
\end{table}

\begin{table}[htb]
\begin{center}
\def\temptablewidth{1.0\textwidth}
{\rule{\temptablewidth}{0.7pt}}
\begin{tabular*}{\temptablewidth}{@{\extracolsep{\fill}}llll}
 %&\multicolumn{3}{c}{L1 scheme }&\multicolumn{3}{c}{Alikhanov's scheme}  \\
% \hline
 L1 scheme &$\sigma \gamma<2-\alpha$& $\sigma \gamma=2-\alpha$ &$\sigma \gamma>2-\alpha$ \\
\hline
\quad Uniform mesh  &$\frac{1-n^{\sigma+\alpha-2}}{2-\alpha -\sigma}$   &    $\ln n $      &$\frac{1-n^{2-\sigma-\alpha}}{\alpha+\sigma-2}  $     \\
\quad Graded mesh  &$\frac{1-n^{\sigma{\gamma}+\alpha-2}}{2-\alpha -\sigma{\gamma}}$     &$\ln n $       &$\ln n $        \\
\quad \Mss{1}  &$\frac{1-(n^*)^{\sigma\gamma+\alpha-2}}{2-\alpha -\sigma\gamma}$     &$\ln n^*$      &$\ln n^*$       \\
\hline
\hline
 Alikhanov's scheme &$\sigma \gamma<3-\alpha$& $\sigma \gamma=3-\alpha$ &$\sigma \gamma>3-\alpha$ \\
\hline
\quad Uniform mesh   &$\frac{1-n^{\sigma+\alpha-3}}{3-\alpha -\sigma}$   &    $\ln n $      &$\frac{1-n^{3-\sigma-\alpha}}{\alpha+\sigma-3}  $        \\
\quad Graded mesh   &$\frac{1-n^{\sigma{\gamma}+\alpha-3}}{3-\alpha -\sigma{\gamma}}$    &$\ln n $       &$\ln n $        \\
\quad \Mss{1}  &$\frac{1-(n^*)^{\sigma\gamma+\alpha-3}}{3-\alpha -\sigma\gamma}$     &$\ln n^*$      &$\ln n^*$       \\
\end{tabular*}
{\rule{\temptablewidth}{0.7pt}}
\end{center}
\caption {$\sigma\in (1,2)$:The formula of factor $\chi_{n,\gamma}^L$ for L1  and $\chi_{n,\gamma}^A$ Alikhanov's schemes with respect to parameters of fractional order $\alpha$, regularity parameter $\sigma$, the grading parameter $\gamma$ and time step $n$.    \label{Table2}}
\end{table}

 The remainder of this paper is organized as follows. In Section \ref{Sec_2}, the main result of $\alpha$-robust error estimate is presented for the general convolution-type scheme in \cite{LMZ19}. After that, the abstract result for general  nonuniform mesh is applied to two typical numerical schemes, i.e., the widely used  L1 scheme and Alikhanov's scheme, and the  $\alpha$-robust  convergence analysis is given with uniform and graded meshes in Sections \ref{Sec_3} and \ref{Sec_4}, respectively.  The paper is ended with conclusion in Section \ref{Sec_5}.

\section{The $\alpha$-robust  error estimate for general discrete schemes\label{Sec_2}}
We here consider a methodology to systematically obtain the $\alpha$-robust error estimate for a discrete convolution-type approximation scheme on nonuniform meshes.
\subsection{The fully discrete general discrete convolution-type  scheme}
Set the generally nonuniform time mesh as $0=t_0<t_1<\cdots<t_N=T$ with the $k$th time-step size~$\tau_k\defeq t_k-t_{k-1}$.  Define $\diff v^k:=v^k-v^{k-1}$,
$t_{n-\theta}\defeq\theta t_{n-1}+(1-\theta)t_n $ and $v^{n-\theta}\defeq\theta v^{n-1}+(1-\theta)v^n$ with an offset parameter $\theta\in[0,1)$.
 Then the convolution-type approximation scheme \cite{LMZ19} of the Caputo derivative is given as
\begin{equation}\label{eq: discrete Caputo}
(\fd{\alpha} v)(t_{n-\theta})\approx \dfd{\alpha}v^{n-\theta}\defeq\sum_{k=1}^n
	A^{(n)}_{n-k}\diff v^k\quad\text{for $1\le n\le N$},
\end{equation}
where   $A^{(n)}_{n-k}$ represents discrete convolution kernels reflecting the convolution structure of the Caputo derivative.
The choice of $\theta =0$ leads to the well-known L1 formula \cite{LLZ18,LX07,OldhamSpanier1974,SOG17,SW06},  and $\theta=\alpha/2$ to the Alikhanov's  approximation scheme~\cite{Alikhanov2015,LMZ21}.

For the spatial discretization, without loss of generality, we take $\Omega = (x_l,x_r)$ for an example. Set $x_i = x_l+ih (0\le i\le M)$ with
  $h\defeq (x_r-x_l)/M$ for integer $M$.  Set a discrete grid
$\overline{\Omega}_h\defeq\{\,x_i\,|\,0\le i\le M\,\}$,
$\Omega_h\defeq\overline{\Omega}_h\cap\Omega$~and
$\partial\Omega_h\defeq\overline{\Omega}_h\cap\partial\Omega$. Denote the space $\mathbb{V}_h$ by the grid functions that vanish on the boundary~$\partial\Omega_h$.
For any $v_h, w_h\in \mathbb{V}_h$, the discrete $L^2$ inner product and the associate $L^2$ norm
are  given as
\begin{align}
\iprod{v_h,w_h}\defeq
h\sum_{x\in\Omega_h}v_h(x)w_h(x), \quad \|v_h\|\defeq\sqrt{\iprod{v_h,v_h}}.
\end{align}
%$$, the $L^2$-norm
% and the $H^1$ semi-norm as
%$$\|v_h\|\defeq\sqrt{\iprod{v_h,v_h}}, \quad \pu{|v_h|_1\defeq\sqrt{h\sum_{x\in\Omega_h}\bra{\nabla_hv_h(x-h/2)}^2}\,,}$$
%where $\nabla_h$ is defined as $\nabla_hv_h(x_{i-\frac12}):=\brab{v_h(x_i)-v_h(x_{i-1})}/h$.
%In the space $\mathbb{V}_h$, the $H^1$ semi-norm $|v_h|_1$ is
%equivalent to the discrete $H^1$-norm $\|v_h\|_1:=\sqrt{\|v_h\|^2+|v_h|_1^2}$, namely
%\[
%C_1 |v_h|_1 \leq  \|v_h\|_1 \leq C_2  |v_h|_1.
%\]
%In this sense, the following estimate of $|v_h|_1$ is also called the $H^1$-norm estimate.
%

We use the second-order approximation  $\Delta_h$ to discretize Laplace operator $\Delta$ and
set  $f_h^n:= f(x_h, t_n)$. The fully discrete scheme for problem \eqref{EQ_equation} is given as
\begin{equation}\label{EQ_GDscheme}
\begin{aligned}
\DCaputo U_h^{n-\theta} &= \Delta_h U_h^{n-\theta}+\kappa U_h^{n-\theta}+f_h^{n-\theta}, && x_h\in \Omega_h,\\
U_h^0 &= u_0(x_h),&& x_h\in \Omega_h,\\
U_h^n & = 0,&& x_h\in \partial \Omega_h, n\geq 0.
\end{aligned}
\end{equation}
The error estimate of scheme \eqref{EQ_GDscheme}
requires the following three assumptions:

\begin{description}
\item[A1.] The discrete kernels are monotone, i.e.,
\[
A^{(n)}_0\ge A^{(n)}_1\ge A^{(n)}_2\ge\cdots A^{(n)}_{n-1}>0
\quad\text{for $1\le n\le N$.}
\]
\item[A2.] There is a constant~$\pi_A>0$ such that the discrete kernels satisfy
a lower bound
\[
A^{(n)}_{n-k}\ge\frac{1}{\pi_A\tau_k}\int_{t_{k-1}}^{t_k}
\omega_{1-\alpha}(t_n-s)\zd{s}\quad\text{for $1\le k\le n\le N$.}
\]
\item[A3.] There is a constant~$\rhomax>0$ such that the step size
ratios~$\rho_k\defeq\tau_k/\tau_{k+1}$ satisfy
\[
\rho_k\le\rhomax\quad\text{for $1\le k\le N-1$.}
\]
\end{description}
\begin{lemma}[\!\cite{LMZ19}]\label{Lemma_theta}
Let assumption \Ass{1} hold and  the parameter $\theta\in [0,1)$. Then it holds for any sequence $\{v^k\}_{k=1}^N$ in $L^2(\Omega)$ that
  \begin{align}
  \sum_{k=1}^nA^{(n)}_{n-k}\nabla_\tau \|v^k\|^2\leq 2\langle \DCaputo v^{n-\theta},  v^{n-\theta}\rangle-d_n(\theta^{(n)}-\theta)\|\DCaputo v^{n-\theta}\|^2, \; 1\leq n\leq N,
  \end{align}
  where $0<d_n<1/A^{(n)}_0$ and $0<\theta^{(n)}<1/2$ are given by
  \begin{align}
    d_n:=\frac{2A^{(n)}_0-A^{(n)}_1}{A^{(n)}_0(A^{(n)}_0-A^{(n)}_1)},\quad \text{and}\quad \theta^{(n)}:=\frac{A^{(n)}_0-A^{(n)}_1}{2A^{(n)}_0-A^{(n)}_1}.
  \end{align}
\end{lemma}
\subsection{The discrete complementary
convolution (DCC) kernels}
We now introduce the definition of DCC kernels $P^{(n)}_{n-k}$, see \cite{LLZ18,LMZ19}, as
    \begin{equation}\label{EQ_PDef1g}
\sum_{j=m}^n P^{(n)}_{n-j}A^{(j)}_{j-m}\equiv1
	\quad\text{for $1\le m\le n\le N$}.
\end{equation}
With this identity, DCC kernel $P^{(n)}_{n-j}$ can be calculated by the following
recursion
\begin{align}\label{discreteConvolutionKernel}
P_{n-k}^{(n)}=\frac{1}{A_{0}^{(k)}}\left\{\begin{array}{cl}
\displaystyle 1\,, &k=n,\vspace{3mm}\\
\displaystyle \sum_{j=k+1}^{n}\brab{A_{j-k-1}^{(j)}-A_{j-k}^{(j)}}P_{n-j}^{(n)}\,, & 1\leq k\leq n-1.
                         \end{array}\right.
\end{align}
It further follows from \eqref{EQ_PDef1g} that
\begin{equation}\label{EQ_PDef2g}\begin{aligned}
\sum_{j=1}^nP_{n-j}^{(n)}\sum_{k=1}^jA_{j-k}^{(j)}\nabla_\tau v^k
=\sum_{k=1}^n\nabla_\tau v^k\sum_{j=k}^nP_{n-j}^{(n)}A_{j-k}^{(j)}
=v^n-v^0, \quad n\geq1.
\end{aligned}
\end{equation}

By taking $v(t)=\omega_{1+\alpha}(t)$ in \cite[formula (2.8) in Lemma 2.1]{LMZ19} such that $\fd{\alpha} v(t)=1$, one can prove that
\begin{equation}
\sum_{j=1}^nP^{(n)}_{n-j}=
    \sum_{j=1}^nP^{(n)}_{n-j}\fd{\alpha}\omega_{1+\alpha}(t_j)
    \le\pi_A\int_0^{t_n}\omega_{1+\alpha}'(s)\zd{s}
    =\pi_A\omega_{1+\alpha}(t_n).
\label{EQ_Pbound}
\end{equation}

%The summation of  DCC kernels is bounded  and satisfies
%\begin{equation}
%\sum_{j=1}^nP^{(n)}_{n-j}\leq \pi_A\omega_{1+\alpha}(t_n),\quad 1\leq n\leq N.\label{EQ_Pbound}
%\end{equation}

A discrete fractional Gr\"{o}nwall inequality in \cite{LMZ19} is presented as follows.
\begin{theorem}[Discrete fractional Gr\"{o}nwall inequality]\label{thm: gronwall}
Assume \Ass{1}--\Ass{3} hold and $0\leq \theta<1$, and the non-negative
sequences  $\{g^n\}_{n=1}^N,$ $\{\lambda_l\}_{l=0}^{N-1}$ and $(v^k)_{k=0}^N$ satisfy
\begin{equation}\label{eq: first Gronwall}
\sum_{k=1}^nA^{(n)}_{n-k}\diff\brab{v^k}^2\le
	\sum_{k=1}^n\lambda_{n-k}\brab{v^{k-\theta}}^2+v^{n-\theta}g^n
	\quad\text{for $1\le n\le N$}.
\end{equation}
If a constant $\Lambda$ satisfies
$\Lambda\ge\sum_{l=0}^{N-1}\lambda_l$ and the maximum step size satisfies
\begin{align}
\max_{1\le n\le N}\tau_n
    \le\frac{1}{\sqrt[\alpha]{2\max(1,\rhomax)\pi_A\Gamma(2-\alpha)\Lambda}}\,, \label{EQ_maxtau}
\end{align}
it holds that
\begin{align}
v^n \leq E_\alpha \big(2\max\{1,\rho\}\pi_A\Lambda t_n^\alpha\big)\Big(v^0+\max_{1\leq k\leq n}\sum_{j=1}^kP^{(k)}_{k-j} g^h\Big),
\end{align}
where $E_\alpha(\cdot):=\sum_{j=0}^\infty \frac{(\cdot)^j}{\Gamma(j\alpha+1)}$ represents the special Mittag-Leffler function.
\end{theorem}

 \subsection{The error estimate of convolution-type scheme  \eqref{EQ_GDscheme}}
Let $u_h^n:=u(x_h, t_n)$ and the error function $e_h^n=u_h^n-U_h^n$ for $x_h\in\bar{\Omega}_h$. Then $e_h^n$ solves
\begin{align}\label{nonuniformL1Scheme-error1}
\bra{\dfd{\alpha}-\Delta_h-\kappa}e_h^{n}
=\mathcal{R}_h^n,\quad x_h\in\Omega_h,\;1\leq n\leq N,
\end{align}
where $\mathcal{R}_h^n := (R_t)_h^n+(R_l)_h^n+(R_s)_h^n$ with
\begin{align}
  (R_t)_h^n&:= \Caputo u^{n-\theta}-\DCaputo u^{n-\theta},\quad  (R_s)_h^n:=\Delta u(x_h,t_{n-\theta})-\Delta_h u(x_h,t_{n-\theta}),\\
  (R_l)_h^n&:=\Delta_h(u (x_h,t_{n-\theta})-u_h^{n-\theta})+\kappa(u (x_h, t_{n-\theta})-u_h^{n-\theta}),\quad x_h\in \Omega_h. \label{EQ_rl}
\end{align}

Taking inner products with $e_h^{n-\theta}$ on both sides of  \eqref{nonuniformL1Scheme-error1}, using Lemma \ref{Lemma_theta} and Theorem \ref{thm: gronwall}, then assuming  $0\leq \theta<\theta^{(n)}$ and maximum step size satisfying \eqref{EQ_maxtau} with $\Lambda = 2\kappa_+$ ($\kappa_+ := \max\{\kappa,0\}$), one has with $C= 2E_\alpha(4\max(1,\rho)\pi_A\kappa t_n^\alpha)$
\begin{align}
\|e_h^n\| \leq &C\Big(\|e_h^0\| +2\max_{1\leq k\leq n}\sum_{j=1}^k P^{(k)}_{k-j}\big(\|(R_s)_h^j\|+\|(R_l)_h^j\|+\|(R_t)_h^j\|\big)\Big).\label{EQ_ger}
\end{align}
The proof of \eqref{EQ_ger} is given in \cite{LMZ19}. The further estimate of \eqref{EQ_ger} requires the choice of the discrete kernels $A^{(n)}_{n-k}$, the discretization of  Laplace operator and  the regularity of  exact solution $u$. For time-fractional diffusion equations, a common regularity assumption in time (cf. \cite{LMZ19,LMZ21,LLZ18}) is given as
\begin{align}
|\partial_t^{l} u^{(l)}(t)|&\leq C_u (1+t^{\sigma-l}), \quad 0<t\leq T, \;l=1,2,\label{EQ_2153}
\end{align}
with the regularity parameter $\sigma\in (0,1)\cup (1,2)$.

It is remarkable that the solution to \eqref{EQ_equation} is smooth at initial time if $\sigma\geq 2$. In this situation,  the optimal convergence order can also be achieved by the same techniques.

Without loss of generality, we assume the solution of \eqref{EQ_equation} in space satisfies
\begin{align}
\mynormb{\partial_{x}^{(4)}u(t)}_{L^\infty}&\leq C_u. \label{SR1}
\end{align}
This spatial regularity can be reduced if the finite element method is used in space \cite{LLZ18}.

We now consider the global consistency error $\sum_{j=1}^k P^{(k)}_{k-j}\big(\|(R_s)_h^j\|+\|(R_l)_h^j\|+\|(R_t)_h^j\|\big)$  in \eqref{EQ_ger} item by item.
\begin{itemize}
\item Under the spacial regularity assumption \eqref{SR1}, the Taylor expansion directly shows that $\|(R_s)^j\|\leq C_uh^2$. Combining with \eqref{EQ_Pbound}, we further have
\begin{equation}
\sum_{j=1}^n P^{(n)}_{n-j}\|(R_s)_h^j\|\leq C_u t_n^\alpha h^2.\label{EQ_1623}
\end{equation}
\item Set $v(t) = (\Delta_h +\kappa)u(x_h,t)$. Using the Taylor formula with integral remainder, we have the integral representation
\begin{equation}
(R_l)^j_h = -\theta\int^{t_{j-\theta}}_{t_{j-1}}(s-t_{j-1})v''(s)\zd s-(1-\theta)\int_{t_{j-\theta}}^{t_j}(t_j-s)v''(s)\zd s, \quad 1\leq j\leq N.
\end{equation}
It is easy to verify that $(R_l)_h^j = 0$ if $\theta=0$, and if $0<\theta<1$,
 we have with the regularity assumption \eqref{EQ_2153} and   appropriate smoothness in space
\begin{equation}
\|(R_l)^j_h\| \leq C_u \tau_1^\sigma/\sigma \quad\text{and}\quad \|(R_l)^j_h\|\leq C_u t_{j-1}^{\sigma-2}\tau_j^2,  \quad 2\leq j\leq N.\label{EQ_2208}
\end{equation}
Following \eqref{discreteConvolutionKernel}, \Ass{2}  that $P^{(n)}_{n-1}\leq 1/A^{(n)}_0\leq \Gamma(2-\alpha)\pi_A\tau_1^\alpha$, and combining with \eqref{EQ_Pbound} and \eqref{EQ_2208}, we have
\begin{align}
\sum_{j=1}^n P^{(n)}_{n-j}\|(R_l)_h^j\|&=P^{(n)}_{n-1}\|(R_l)_h^1\|+\sum_{j=2}^n P^{(n)}_{n-j}\|(R_l)_h^j\|\nonumber\\
&\leq C_u \big(\tau_1^{\sigma+\alpha}/\sigma+t_n^\alpha\max_{2\leq k\leq n}t_{k-1}^{\sigma-2}\tau_k^2\big).\label{EQ_1619}
\end{align}
\item The estimate of $\sum_{j=1}^k P^{(k)}_{k-j}\|(R_t)_h^j\|$ is the most complicate comparing with the above two terms.  The reasons have two folds: (i) Noting that $(R_t)_h^j$ represents the local truncation error of scheme \eqref{eq: discrete Caputo} for Caputo derivative, it implies that the estimate depends on the total error contributions
from discretization errors over the $k$~current intervals $[t_{j-1},t_{j}]$
($1\leq j\leq k$). Thus one needs to calculate the discretization
errors over the $\frac{1}{2}k(k-1)$ historic intervals
$[t_{0},t_{j-1}]$ ($2\leq j\leq k$). Hence it is very cumbersome to calculate the consistence error, also is the reason why it is challenging to deal with general nonuniform meshes. (ii) Because the explicit formula of $P_{n-j}^{(n)}$ is unknown, it brings another challenge to estimate on generally nonuniform meshes.

To circumvent this difficulty,  the concept of
\emph{error convolution structure} (ECS) proposed in \cite{LMZ21,LLZ18} is introduced to simplify the estimate of global consistency error.  Specifically, the ECS tells us that  $(R_t)_h^k$  be bounded by the following discrete convolution structure
\begin{equation}
  \text{({\bf ESC} hypothesis)} \quad \|(R_t)_h^k\|\leq A^{(k)}_0\Gloc^k+\sum_{j=1}^{k-1}(A^{(k)}_{k-j-1}-A^{(k)}_{k-j})\Ghis^j.  \label{EQ_ESC}
\end{equation}
Here $\Gloc^k$ and $\Ghis^j$ respectively represent the interpolation errors on the local  subinterval $[t_{k-1}, t_{k-\theta}]$ and history $[t_0,t_{k-1}]$, which depend on the specific approximation schemes and will be presented in details later.
With the {\bf ESC} hypothesis \eqref{EQ_ESC}, the global consistency error can be reduced to
\begin{align}
\sum_{k=1}^n P^{(n)}_{n-k}\|(R_t)_h^k\|&\leq \sum_{k=1}^nP^{(n)}_{n-k}A_0^{(k)}\Gloc^k+\sum_{k=1}^{n}P^{(n)}_{n-k}
\sum_{j=1}^{k-1}(A^{(k)}_{k-j-1}-A^{(k)}_{k-j})\Ghis^j\nonumber\\
&\leq \sum_{k=1}^nP^{(n)}_{n-k}A_0^{(k)}\Gloc^k+\sum_{j=1}^{n-1}\Ghis^j
\sum_{j=1}^{k-1}(A^{(k)}_{k-j-1}-A^{(k)}_{k-j})P^{(n)}_{n-k}\nonumber\\
&=\sum_{k=1}^nP^{(n)}_{n-k}A_0^{(k)}\Gloc^k+\sum_{k=1}^{n-1}P^{(n)}_{n-k}A_0^{(k)}\Ghis^k \nonumber\\
&\leq\sum_{k=1}^nP^{(n)}_{n-k}A_0^{(k)}(\Gloc^k+\Ghis^k) := \sum_{k=1}^nP^{(n)}_{n-k} \mathcal{G}^k
, \label{EQ_2045}
\end{align}
where \eqref{discreteConvolutionKernel} is used in the first identity.  The bound \eqref{EQ_2045} shows that the ECS technique can significantly reduce evaluation of the global consistency error from  the $\frac{1}{2}n(n-1)$ historic intervals to $n$ intervals.  Hence, ECS is a useful tool for the convergence analysis. We now consider the $\alpha$-robust  error estimate for \eqref{EQ_2045} by introducing the following lemma.
\begin{lemma}\label{Lemma_npe}
Let  $\{A^{(n)}_{n-k}\}_{k=1}^n$ and $\{P^{(n)}_{n-k}\}_{k=1}^n$ be defined by
 \eqref{eq: discrete Caputo}  and \eqref{EQ_PDef1g} respectively. If \Ass{1} and \Ass{2} hold, for any positive sequence $\{\nu_k\}_{k=1}^n$, one has
\begin{align}
\sum_{k=1}^{n}P_{n-k}^{(n)}  \nu_k
\leq \Gamma(2-\alpha)  \pi_A\sum_{j=1}^n\tau_j \max_{j\leq k\leq n}t_k^{\alpha-1} \nu_k.\label{EQ_npe}
\end{align}
\end{lemma}
\begin{proof}
It follows from \Ass{2} that
\begin{equation*}
\sum_{j=1}^kA^{(k)}_{k-j}\tau_j\geq \frac1{\pi_A}\int_{0}^{t_k}\omega_{1-\alpha}(t_k-s)\zd s =\frac1{\pi_A} \omega_{2-\alpha}(t_k)=\frac{t_k^{1-\alpha}}{\pi_A\Gamma(2-\alpha) }.
\end{equation*}
Thus, we have
\begin{align*}
\sum_{k=1}^{n}P_{n-k}^{(n)} \nu_k \cdot 1
&\leq  \sum_{k=1}^{n}P_{n-k}^{(n)}\nu_k
\cdot \Gamma(2-\alpha) \pi_A t_k^{\alpha-1}\sum_{j=1}^kA^{(k)}_{k-j}\tau_j \\
&\leq \Gamma(2-\alpha)   \pi_A\sum_{j=1}^n\tau_j\sum_{k=j}^nP_{n-k}^{(n)}A^{(k)}_{k-j}t_k^{\alpha-1} \nu_k\\
&\leq \Gamma(2-\alpha)  \pi_A \sum_{j=1}^n\tau_j \Big(\max_{j\leq k\leq n}t_k^{\alpha-1} \nu_k\Big) \underbrace{\sum_{k=j}^nP_{n-k}^{(n)}A^{(k)}_{k-j}}_{=1},
\end{align*}
where the last inequality uses the identity \eqref{EQ_PDef1g}. The proof is completed.
\end{proof}
\medskip

Applying Lemma \ref{Lemma_npe} by taking $v_k = \mathcal{G}^k$ in \eqref{EQ_npe} to \eqref{EQ_2045}, one has
\begin{align}
&\sum_{j=1}^k P^{(k)}_{k-j}\|(R_t)_h^j\|\nonumber\\
\leq& \Gamma(2-\alpha)  \pi_A\sum_{j=1}^n\tau_j \max_{j\leq k\leq n}t_k^{\alpha-1} \mathcal{G}^k \nonumber\\
=&\Gamma(2-\alpha)  \pi_A\Big(\sum_{j=2}^n\tau_j \max_{j\leq k\leq n}t_k^{\alpha-1}\mathcal{G}^k+\max\big\{\tau_1\max_{2\leq k\leq n}t_k^{\alpha-1}\mathcal{G}^k, \tau_1^\alpha \mathcal{G}^1\big\}\Big)\nonumber\\
\leq &\Gamma(2-\alpha)  \pi_A\Big(\sum_{j=2}^n\tau_j \max_{j\leq k\leq n}t_k^{\alpha-1}\mathcal{G}^k+ \rho\tau_2\max_{2\leq k\leq n}t_k^{\alpha-1}\mathcal{G}^k+\tau_1^\alpha \mathcal{G}^1\Big)\nonumber\\
\leq &\Gamma(2-\alpha)  \pi_A\Big((1+\rho)\sum_{j=2}^n\tau_j \max_{j\leq k\leq n}t_k^{\alpha-1}\mathcal{G}^k+ \tau_1^\alpha \mathcal{G}^1\Big)\nonumber\\
\leq &(1+\rho)\Gamma(2-\alpha)  \pi_A\Big(\sum_{j=2}^n\tau_j \max_{j\leq k\leq n}t_k^{\alpha-1}\mathcal{G}^k+ \tau_1^\alpha \mathcal{G}^1\Big) :=  (\rho+1)\Xi^n_t.
 \label{RtE}
\end{align}
\end{itemize}

%\subsubsection{Two examples for error estimate \eqref{EQ_ger2}}

%\begin{remark}
%Note that  there is a factor $1/(1-\alpha)$  in error estimates \eqref{EQ_esL1} and \eqref{EQ_esAl} and  the factor will blow up as $\alpha\to 1^-$. It implies  the estimates \eqref{EQ_esL1} and \eqref{EQ_esAl} are $\alpha$-nonrobust.  Actually,  the factor $1/(1-\alpha)$ in these estimates  arises from a temporal maximum-norm   bound (see \cite[(3.6)]{LLZ18} and
% \cite[Lemma 3.4]{LMZ21}) of the global consistency error $2\max_{1\leq k\leq n}\sum_{j=1}^k P^{(k)}_{k-j}A_0^{(j)}(\Gloc^j+\Ghis^j)$ in \eqref{EQ_ger2}. To obtain the $\alpha$-robust error estimate, it is necessary to revisit the bound of the global consistency error.
%\end{remark}
By combining \eqref{EQ_ger}, \eqref{EQ_1623}, \eqref{EQ_1619} and \eqref{RtE}, we have our main results of $\alpha$-robust  error estimate for general scheme \eqref{EQ_GDscheme} as follows.
\begin{theorem}\label{Theorem_main_result}
 Assume   \Ass{1}--\Ass{3},  the regularity assumptions \eqref{EQ_2153} and \eqref{SR1} hold. If the maximum step size satisfies
 \begin{align}
\max_{1\le n\le N}\tau_n
    \le\frac{1}{\sqrt[\alpha]{4\max(1,\rhomax)\pi_A\Gamma(2-\alpha)\kappa_+}},
\end{align}
 then it holds with  $0<\theta<\theta^{(n)}$  and $C=C_uE_\alpha(4\max(1,\rho)\pi_A\kappa t_n^\alpha)$ that
\begin{equation}
\|e_h^n\|\leq C\Big(\|e_h^0\|+\frac{\tau_1^{\sigma+\alpha}}{\sigma}+t_n^\alpha h^2+t_n^\alpha\max_{2\leq k\leq n}t_{k-1}^{\sigma-2}\tau_k^2
+(\rho+1)\max_{1\leq k\leq n}\Xi^k_t\Big),\label{EQ_ger3}
\end{equation}
and it holds with $\theta =0$ that
\begin{equation}
\|e_h^n\|\leq C\Big(\|e_h^0\|+t_n^\alpha h^2+(\rho+1)\max_{1\leq k\leq n}\Xi^k_t\Big),\label{EQ_ger31}
\end{equation}
where $\Xi^k_t$ is defined in \eqref{RtE} by
\begin{equation} \label{EQ_Xit}
\Xi^k_t= \Gamma(2-\alpha)  \pi_A\Big(\sum_{j=2}^n\tau_j \max_{j\leq k\leq n}t_k^{\alpha-1}\mathcal{G}^k+ \tau_1^\alpha \mathcal{G}^1\Big).
\end{equation}
\end{theorem}

\begin{remark}
We point out that there exists a factor $1/(1-\alpha)$ for the error estimates of \eqref{EQ_2045} in previous theory \cite{LLZ18,LMZ21,LMZ19}. But, our new estimate \eqref{EQ_Xit} for \eqref{EQ_2045} is bounded by $\Gamma(2-\alpha)$ instead of $1/(1-\alpha)$.   As $\alpha\to 1^-$, one has $\Gamma(2-\alpha)\rightarrow 1$ instead of $1/(1-\alpha) \rightarrow \infty$. Thus, our new estimate circumvents the factor blow-up phenomenon. Next, we need further study if the remainder term in \eqref{EQ_Xit} exists the factor blow-up phenomenon. Specifically, we take L1 and Alikhanov's schemes for examples to revisit the main theorem \eqref{Theorem_main_result}, and present the detailed factors how to depend on $\alpha$, $\sigma$ and mesh parameters in the convergence analysis by estimating \eqref{EQ_Xit} with various meshes.
\end{remark}

\section{Application to L1 scheme\label{Sec_3}}
We first consider the $\alpha$-robust error estimate to L1 scheme \cite[p.~140]{OldhamSpanier1974} with  the discrete convolution kernels given as
\begin{equation} \label{L1s}
A^{(n)}_{n-k}\defeq\frac{1}{\tau_k}\int_{t_{k-1}}^{t_k}
    \omega_{1-\alpha}(t_n-s)\zd{s}.
\end{equation}
As proved in \cite{LLZ18}, one can verify the above assumptions hold as follows:
\begin{enumerate}
  \item Assumptions \Ass{1} holds and \Ass{2} holds with $\pi_A = 1$, and $\rho>0$ in \Ass{3};
  \item $\theta = 0$;
  \item The local and history interpolation errors $\Gloc^k$ and $\Ghis^k$ are in the form of
\begin{equation} \label{DG}
\Ghis^k\triangleq \int_{t_{k-1}}^{t_k}\bra{s-t_{k-1}}\|u_{tt}(s)\|\zd s\,,
\; 1\leq k\leq n-1 , \; \text{and} \quad \Gloc^n = \Ghis^n;
\end{equation}
\item Using the integral mean-value theorem, it is easy to verify
\begin{equation} \label{MVT}
A_{n-k+1}^{(n)}<\omega_{1-\alpha}(t_n-t_{k-1})<A_{n-k}^{(n)}.
\end{equation}
\end{enumerate}
Under the regularity condition \eqref{EQ_2153} and $u\in C^2((0,T])$, it follows from \eqref{DG} that
\begin{align}
 \Gloc^1\leq C_u\tau_1^{\sigma}/\sigma, \quad \text{and} \quad \Gloc^k\leq C_ut_{k-1}^{\sigma-2}\tau_k^2 \quad \text{for}
\quad 2\leq k\leq n. \label{EQ_Gke}\end{align}

Combining with \eqref{EQ_Gke} and \eqref{EQ_Xit}, $\Xi^k_t$ can be bounded by
\begin{align}
\Xi^n_t&\leq C_u\Gamma(2-\alpha)\Big(\frac{\tau_1^\sigma}{\sigma}+\sum_{j=2}^n\tau_j\max_{j\leq k\leq n}(t_{k-1}^{\sigma-2}t_k^{\alpha-1}\tau_k^{2-\alpha})\Big).\label{EQ_alpharobust}
\end{align}
\begin{remark} \label{Remark_L1}
To make our analysis clear how to circumvent the $\alpha$ blow-up phenomenon, we now recall the techniques used in \cite{LLZ18}, which follows from \eqref{EQ_2045} that
\begin{align} \label{P1}
\sum_{j=1}^nP_{n-j}^{(n)}\absb{(R_t)^j}\leq\, 2P_{n-1}^{(n)}A_{0}^{(1)}\Gloc^1+2\sum_{k=2}^{n}P_{n-k}^{(n)}A_{k-2}^{(k)}\frac{A_{0}^{(k)}}{A_{k-2}^{(k)}}\Gloc^k.
\end{align}
Using the relationship in \eqref{MVT}, it arrives at
\begin{align*}
\frac{A_{0}^{(k)}}{A_{k-2}^{(k)}}<\frac{\omega_{2-\alpha}(\tau_k)}{\tau_k\,\omega_{1-\alpha}(t_k-t_1)}
=\frac{(t_k-t_1)^{\alpha}}{1-\alpha}\tau_k^{-\alpha},\quad 2\leq k\leq n\,.
\end{align*}
Applying the bounds in \eqref{EQ_Gke}, taking the maximum value of $\frac{A_{0}^{(k)}}{A_{k-2}^{(k)}}\Gloc^k$ for $k = 2,\cdots,n$ in \eqref{P1}, using  facts in \eqref{EQ_PDef1g} that
$P_{n-1}^{(n)}A_{0}^{(1)}\leq1$ and $\sum_{k=2}^{n}P_{n-k}^{(n)}A_{k-2}^{(k)}=1$, one has
\begin{align}\label{globalL1ConsistenceError1}
\sum_{j=1}^nP_{n-j}^{(n)}\abs{(R_t)^j} \leq C_u\braB{\frac{1}{\sigma}\tau_1^{\sigma}+\frac{1}{1-\alpha}\max_{2\leq k\leq n}(t_{k}-t_1)^{\alpha}t_{k-1}^{\sigma-2}\tau_k^{2-\alpha}},\quad n\geq1,
\end{align}
which hence introduces the factor of $1/(1-\alpha)$. But our bound \eqref{EQ_Gke} avoids the factor.
\end{remark}

We next investigate the  global consistency error with various meshes. Beforehand, we introduce  two factors which will appear  frequently in the error bounds. Set \begin{equation}
\beta_\gamma = \min\{\sigma\gamma,2-\alpha\},\quad \gamma\geq 1.\label{EQ_betagamma}
\end{equation}
Then we denote by
\begin{equation}
  \varsigma_{n,\gamma}^L:= n^{\alpha-2+\beta_\gamma}\int_{1/n}^{1} s^{\sigma\gamma+\alpha-3} \zd s= \left\{\begin{array}{cc}
\frac{1-n^{\sigma\gamma+\alpha-2}}{2-\alpha-\sigma\gamma},& \sigma\gamma< 2-\alpha,\\
\ln n,&\sigma\gamma=2-\alpha,\\
\frac{1-n^{2-\alpha-\sigma\gamma}}{\sigma\gamma+\alpha-2},&\sigma\gamma>2-\alpha,
\end{array}\right.\label{EQ_etabar}
\end{equation}
and
\begin{align}
\zeta_{n,\gamma}^L :=n^{\beta_\gamma+\alpha-2}\int_{1/n}^1s^{\beta_\gamma+\alpha-3}\zd s=\left\{\begin{array}{cc}
\frac{1-n^{\sigma\gamma+\alpha-2}}{(2-\alpha)-\sigma\gamma},&\sigma\gamma<2-\alpha,\\
\ln n,&\sigma\gamma\geq 2-\alpha.
\end{array}\right.\label{EQ_1708}
\end{align}
For convenience, we further define
\begin{equation}\label{EQ_cl}
  \chi_{n,\gamma}^L:= \left\{\begin{array}{cc}
\varsigma_{n,\gamma}^L,&0<\sigma<1,\\
\zeta_{n,\gamma}^L,&1<\sigma<2.
\end{array}\right.
\end{equation}
\begin{remark}\label{EQ_ramark1} It follows from \eqref{EQ_etabar} and \eqref{EQ_1708}  that $\varsigma_{n,\gamma}^L, \zeta_{n,\gamma}^L$  are continuous with respect to $\alpha, \sigma$, and $\varsigma_{n,\gamma}^L, \zeta_{n,\gamma}^L \to \ln n$ as $\sigma+\alpha\to 2$, which implies $\varsigma_{n,\gamma}^L$ and $\zeta_{n,\gamma}^L$ would not blow up, and can be uniformly bounded by $\ln n$ for any $\alpha$ and $\sigma$. This is, it holds
      \begin{align*}
        \varsigma_{n,\gamma}^L &=\zeta_{n,\gamma}^L= n^{\sigma\gamma+\alpha-2}\int_{1/n}^{1} s^{\sigma\gamma+\alpha-3} \zd s\leq \int_{1/n}^{1} s^{-1} \zd s = \ln n, && \text{as } \sigma\gamma<2-\alpha, \\
        \varsigma_{n,\gamma}^L &= \int_{1/n}^{1} s^{\sigma\gamma+\alpha-3} \zd s\leq \int_{1/n}^{1} s^{-1} \zd s = \ln n,  &&\text{as } \sigma\gamma>2-\alpha.
      \end{align*}
%where one uses the facts that $\sigma<2, \alpha<1$ and $t_n^{\sigma+\alpha-2}\leq \max\{1,T^{\sigma+\alpha-2}\}\leq T+1 $.
\end{remark}
\subsection{The global consistency error $\Xi^n_t$ on uniform mesh}
We first consider uniform mesh by taking $\tau = T/N$ and $t_k= k\tau$, thus estimate \eqref{EQ_alpharobust} has
\begin{equation}\label{EQ_1751}
\Xi^n_t\leq C_u(\frac1{\sigma}\tau^\sigma+\sum_{j=2}^n\tau^{3-\alpha}t_{j-1}^{\sigma+\alpha-3}).
\end{equation}
The second term in right hand of \eqref{EQ_1751} can be estimated by
\begin{align}
\sum_{j=2}^n\tau^{3-\alpha}t_{j-1}^{\sigma+\alpha-3}&\leq \tau^\sigma+\tau^{2-\alpha}\sum_{j=2}^{n}\tau t_j^{\sigma+\alpha-3} \nonumber \\ &=\tau^{\sigma}\Big(1+n^{\sigma+\alpha-2}\sum_{j=2}^{n}
\frac1{n}\big(\frac{j}{n}\big)^{\sigma+\alpha-3}\Big)\nonumber
\\
&\leq \tau^\sigma\Big(1+n^{\sigma+\alpha-2}\int_{1/n}^{1} s^{\sigma+\alpha-3} \zd s\Big)\; =\tau^\sigma+\varsigma_{n,1}^Lt_n^{\sigma-\beta_1}\tau^{\beta_1},\label{EQ_1750}
\end{align}
where $\beta_1$ and  $\varsigma_{n,1}^L$ are defined in \eqref{EQ_betagamma} and \eqref{EQ_etabar}, repectively.
%the factor $\bar{\varsigma}_n^*$ is given as
%\begin{equation}
%  \bar{\varsigma}_n^*:= n^{\alpha-2+\beta_1}\int_{1/n}^{1} s^{\sigma+\alpha-3} \zd s= \left\{\begin{array}{cc}
%\frac{1-n^{\sigma+\alpha-2}}{2-\alpha-\sigma},& \sigma< 2-\alpha\\
%\ln n,&\sigma=2-\alpha\\
%\frac{1-n^{2-\sigma-\alpha}}{\alpha+\sigma-2},&\sigma>2-\alpha.
%\end{array}\right.
%\end{equation}
Inserting the estimate \eqref{EQ_1750} into \eqref{EQ_1751}, we arrive at
\begin{equation}
\Xi^n_t\leq C_u\Big((\frac1{\sigma}+1)\tau^\sigma+\varsigma_{n,1}^Lt_n^{\sigma-\beta_1}\tau^{\beta_1}\Big).\label{EQ_L1uni}
\end{equation}

\begin{remark}
For the case of $\kappa\leq 0$ and $\sigma=\alpha$, an  $\alpha$-robust bound of the  global consistency error on uniform mesh has been presented in \cite{CS21} as
\begin{align}
\sum_{j=1}^k P^{(k)}_{k-j}\|(R_t)_h^j\|\leq CK_{2-\alpha,n}N^{-1}t_n^{\alpha-1},
\end{align}
where  $$K_{\beta,n}:=\Big\{\begin{array}{cc}
1+\frac{1-n^{1-\beta}}{\beta-1} & \text{if}\; \beta\neq 1,\\
1+\ln n&\text{if}\; \beta=1.
\end{array}.$$
It is easy to check $K_{2-\alpha,n}$ is  almost the same as $ \varsigma_{n,1}^L$ for the special case of   $\sigma=\alpha$. It implies that our factor $\varsigma_{n,1}^L$ is consistent with  $K_{2-\alpha,n}$ .
\end{remark}

\subsection{The global approximation error $\Xi^n_t$ on graded meshes}
%A typical nonuniform mesh is the  graded mesh defined as $t_k=T\bra{k/N}^{\gamma}$ with $\gamma \geq 1$. If $\gamma=1$, the mesh will be uniform. 
We now consider the global approximation error $\Xi^n_t$ on the graded mesh.
It is easy to verify $\tau_k\leq \gamma Tk^{\gamma-1}N^{-\gamma}$.

Note that $\sigma>0$, it follows from \eqref{EQ_alpharobust} that
\begin{align}
\Xi^n_t&\leq C_uT^\sigma\Big(\frac{N^{-\gamma\sigma}}{\sigma}+\frac{\gamma^{2-\alpha}}{T}\sum_{j=2}^n\tau_j\max_{j\leq k\leq n}(\frac{k}{N})^{-\gamma}(\frac{k}{k-1})^{(2-\sigma)\gamma}k^{\sigma\gamma-(2-\alpha)}N^{-\sigma\gamma}\Big)\nonumber\\
&\leq C_uT^\sigma\Big(\frac{N^{-\gamma\sigma}}{\sigma}+4^\gamma\gamma^{2-\alpha} T^{-1}\sum_{j=2}^n\tau_j\max_{j\leq k\leq n}(\frac{k}{N})^{-\gamma}k^{\sigma\gamma-(2-\alpha)}N^{-\sigma\gamma}\Big).
 \label{EQ_1424}
\end{align}
Note that $(\sigma-1)\gamma-(2-\alpha)<0$ for any $\sigma,\alpha\in (0,1)$. Thus, for $\sigma\in (0,1)$, \eqref{EQ_1424} can be reduced to
\begin{equation}
\Xi^n_t\leq C_uT^\sigma N^{-\sigma\gamma}\Big(\frac{1}{\sigma}+4^\gamma\gamma^{3-\alpha}\sum_{j=2}^nj^{\sigma\gamma-(3-\alpha)}\Big).\label{EQ_1532}
\end{equation}
For any $\xi$, the direct calculation shows that
\begin{align}
\sum_{j=2}^{n-1}j^{-\xi-1} %&= n^{-\xi}\sum_{j=2}^n\frac1{n}\left(\frac{j}{n}\right)^{-\xi-1}
\leq n^{-\xi}\int_{1/n}^1s^{-\xi-1}\zd s= \left\{\begin{array}{cc}
\frac{1}{\xi}(1-n^{-\xi}),&\xi\neq0\\
\ln n,&\xi= 0.
\end{array}\right.\label{EQ_0116}
\end{align}
Then setting $\xi = 2-\alpha-\sigma\gamma$ in \eqref{EQ_0116} and inserting the resulting into \eqref{EQ_1532}, one has
\begin{align}
\Xi^n_t&\leq C_u\Big(\frac{T^\sigma}{\sigma}N^{-\sigma\gamma}+4^\gamma\gamma^{3-\alpha}
T^{\frac{\beta_\gamma}{\gamma}}t_n^{\sigma-\frac{\beta_\gamma}{\gamma}}(n^{\beta_\gamma-(3-\alpha)}+
\varsigma_{n,\gamma}^L)N^{-\beta_\gamma}\Big)\nonumber\\
&\leq C_uC_\gamma T^\sigma\Big(\frac{1}{\sigma}+
\varsigma_{n,\gamma}^L\Big)N^{-\min\{\gamma\sigma,2-\alpha\}},\label{EQ_XiL1g1}
\end{align}
where $C_\gamma$ is a constant and $\beta_\gamma, \varsigma_{n,\gamma}^L$ are defined in \eqref{EQ_betagamma} and  \eqref{EQ_etabar}, respectively.

On the other hand, for $\sigma\in (1,2)$, it follows from \eqref{EQ_1424} that
\begin{align}
\Xi^n_t&\leq C_uT^\sigma\Big(\frac{N^{-\gamma\sigma}}{\sigma}+\frac{\gamma^{2-\alpha}}{T}\sum_{j=2}^n\tau_j\max_{j\leq k\leq n}(\frac{k}{N})^{-\gamma}(\frac{k}{k-1})^{(2-\sigma)\gamma}k^{\sigma\gamma-(2-\alpha)}N^{-\sigma\gamma}\Big)\nonumber\\
&\leq C_uT^\sigma\Big(\frac{N^{-\gamma\sigma}}{\sigma}+4^\gamma\gamma^{2-\alpha} T^{-1}\sum_{j=2}^n\tau_j\max_{j\leq k\leq n}(\frac{k}{N})^{-\gamma}k^{\sigma\gamma-(2-\alpha)}N^{-\sigma\gamma}\Big)\nonumber\\
&\leq C_uT^\sigma\Big(\frac{N^{-\gamma\sigma}}{\sigma}+4^\gamma\gamma^{3}N^{-\beta_1}\sum_{j=2}^nj^{-1}
\max_{j\leq k\leq n}(\frac{k}{N})^{\sigma\gamma-\beta_\gamma}k^{\beta_\gamma-(2-\alpha)}\Big)\nonumber\\
&\leq C_uT^\sigma\Big(\frac{N^{-\gamma\sigma}}{\sigma}+4^\gamma\gamma^{3}N^{-\beta_\gamma}\sum_{j=2}^nj^{\beta_\gamma-(3-\alpha)}
\max_{j\leq k\leq n}(\frac{k}{N})^{\sigma\gamma-\beta_\gamma}\Big)\nonumber \\
&\leq C_uT^\sigma\Big(\frac{N^{-\gamma\sigma}}{\sigma}+4^\gamma\gamma^{3}N^{-\beta_\gamma}\sum_{j=2}^nj^{\beta_\gamma-(3-\alpha)}
\Big) .
 \label{EQ_1803}
\end{align}
%$$\varsigma_n = \sum_{j=2}^nj^{-1}\max_{j\leq k\leq n}(\frac{k}{N})^{\sigma\gamma-\min(\sigma\gamma, 2-\alpha)}k^{\min(\sigma\gamma, 2-\alpha)-(2-\alpha)}.$$
For any $\xi$, the direct calculation shows that
\begin{align}
\sum_{j=2}^{n}j^{-\xi-1} %&= n^{-\xi}\sum_{j=2}^n\frac1{n}\left(\frac{j}{n}\right)^{-\xi-1}
\leq n^{-\xi}\int_{1/n}^1s^{-\xi-1}\zd s= \left\{\begin{array}{cc}
\frac{1}{\xi}(1-n^{-\xi}),&\xi>0\\
\ln n,&\xi= 0.
\end{array}\right.\label{EQ_0117}
\end{align}
Then setting $\xi = 2-\alpha-\min\{\sigma\gamma, 2-\alpha\}$ in \eqref{EQ_0117} and inserting the resulting into  \eqref{EQ_1803}, we have
\begin{align}
\Xi^n_t
\leq&
C_uT^{\sigma} \braB{\frac{1}{\sigma} N^{-\gamma\sigma}+4^\gamma\gamma^{3}\zeta_{n,\gamma}^L N^{-\min\{\gamma\sigma,2-\alpha\}}}\nonumber\\
\leq& C_uC_\gamma T^{\sigma} \braB{\frac{1}{\sigma} +\zeta_{n,\gamma}^L }N^{-\min\{\gamma\sigma,2-\alpha\}},
\label{EQ_XiL1g2}
\end{align}
where $C_\gamma$ is a constant and $\zeta_{n,\gamma}^L$ is defined in \eqref{EQ_1708}.
%Here the factor $\varsigma_n^*$ is given by
%\begin{align}
% \varsigma_n^* :=n^{\beta_\gamma+\alpha-2}\int_{1/n}^1s^{\beta_\gamma+\alpha-3}\zd s=\left\{\begin{array}{cc}
%\frac{1-n^{\sigma\gamma+\alpha-2}}{(2-\alpha)-\sigma\gamma},&\gamma<\frac{2-\alpha}{\sigma},\\
%\ln n,&\gamma\geq\frac{2-\alpha}{\sigma}.
%\end{array}\right.
%\end{align}
Thus,  from \eqref{EQ_XiL1g1}, \eqref{EQ_XiL1g2}, the global consistency error $\Xi^n_t$ on graded meshes can be organized as
\begin{equation}\label{gradedconsistency}
\begin{aligned}
\Xi^n_t
\leq&  C_uC_\gamma T^{\sigma} \braB{\frac{1}{\sigma} +\chi_{n,\gamma}^L }N^{-\min\{\gamma\sigma,2-\alpha\}}, \\
\end{aligned}
\end{equation}
where the factor $\chi_{n,\gamma}^L$ is defined in \eqref{EQ_cl}.
To obtain the optimal convergence, the optimal choice of $\gamma$ is to take
$\gamma_{\texttt{opt}}=\max\{1,(2-\alpha)/\sigma\}.$
%\begin{remark}
%$\varsigma_n^*$ is also continuous with respect to $\alpha, \sigma$ and has an upper bound
%      \begin{align*}
%        \varsigma_n^* = n^{\sigma+\alpha-2}\int_{1/n}^1s^{\sigma+\alpha-3}\zd s\leq \int_{1/n}^{1} s^{-1} \zd s \leq \ln n, \quad \text{for } \sigma<2-\alpha.
%      \end{align*}
%\end{remark}
\begin{remark}
 For the case of $\kappa\leq 0, \sigma=\alpha$, a similar $\alpha$-robust bound of the   global consistency error on graded mesh has been presented in \cite{CS21} as
\begin{align}
\sum_{j=1}^k P^{(k)}_{k-j}\|(R_t)_h^j\|\leq \frac{Ce^{\gamma}\Gamma(1+\ell_N-m^*/\gamma) T^\alpha}{\Gamma(1+\ell_N+\alpha-m^*/\gamma)\Gamma(2-\alpha)}
\left(\frac{t_n}{ T}\right)^{\ell_N+\alpha-m^*/\gamma}N^{-m^*},
\end{align}
where $\ell_N:= 1/(\ln N)$ and $m^*:=\min\{\alpha\gamma, 2-\alpha \}$. This error bound suffers from a minor restriction $\gamma\leq 2(2-\alpha)/\alpha$, but is circumvented in our consistency error \eqref{gradedconsistency}.
\end{remark}

\subsection{The global error $\Xi^n_t$ on general nonuniform mesh}
We now consider the global error $\Xi^n_t$ on general nonuniform mesh \Mss{1}.
%presented in \cite{MM15,B85}  as
%\begin{enumerate}
%\item[\Mss{1}.] There is a constant $C_1,C_{\gamma}$ such that \cy{$\tau_1\geq C_1 \tau^\gamma$} and
%$\tau_k\le C_{\gamma}\tau\min\{1,t_k^{1-1/\gamma}\}$
%for~$1\le k\le N$,  $t_{k}\leq C_{\gamma}t_{k-1}$ and $\tau_k \leq  C_{\gamma} \tau_{k-1}$ for~$2\le k\le N$.
%\end{enumerate}
%\cite{MM15,B85}One can refer to this nonuniform mesh \Mss{1} for solving parabolic equations with nonsmooth initial values \cite{LM21,LM22,LMW21}.
%With the general nonuniform mesh, we now present the $\alpha$-robust estimate of the global consistency error.

With the mesh assumption \Mss{1}, \eqref{EQ_alpharobust} can be reduced to
\begin{align}
\Xi^n_t&\leq C_uC_\gamma\Big(\frac{\tau^{\sigma\gamma}}{\sigma}+\tau^{2-\alpha}\sum_{j=2}^n\tau_j\max_{j\leq k\leq n}t_k^{\sigma-\frac{2-\alpha}{\gamma}-1}\Big).\label{EQ_1648}
\end{align}
It is easy to check that  $\sigma-\frac{2-\alpha}{\gamma}-1<0$ for any $\gamma\geq 1$ if $\sigma\in(0,1)$. Then, \eqref{EQ_1648} can be further reduced to
\begin{align}
\Xi^n_t&\leq C_uC_\gamma\Big(\frac{\tau^{\sigma\gamma}}{\sigma}+\tau^{2-\alpha}\sum_{j=2}^n\tau_jt_j^{\sigma-\frac{2-\alpha}
{\gamma}-1}\Big).\label{EQ_0209}\end{align}
For any $\xi$, the direct calculation shows that
\begin{align}
\sum_{j=2}^{n-1}\tau^{j} t_j^{-\xi-1}
\leq (\rho+1) \int_{\tau_1}^1s^{-\xi-1}\zd s= \left\{\begin{array}{cc}
\frac{\rho+1}{\xi}(\tau_1^{-\xi}-{t_n}^{-\xi}),&\xi>0,\\
(\rho+1)\ln \frac{t_n}{\tau_1},&\xi= 0.
\end{array}\right.\label{EQ_0201}
\end{align}
Setting $\xi = (2-\alpha)/\gamma-\sigma$ in \eqref{EQ_0201} and inserting the resulting into  \eqref{EQ_0209}, we have
\begin{align}
\Xi^n_t&\leq C_uC_\gamma\Big(\frac{\tau^{\sigma\gamma}}{\sigma}+\tau^{2-\alpha}\tau_n t_n^{\sigma-\frac{2-\alpha}{\gamma}-1}+
\tau^{2-\alpha}\sum_{j=2}^{n-1}\tau_jt_j^{\sigma-\frac{2-\alpha}
{\gamma}-1}\Big)\nonumber\\
&\leq C_uC_\gamma\Big(\frac{\tau^{\sigma\gamma}}{\sigma}+\tau^{\beta_\gamma}\big(\frac{\tau_1}{t_n}\big)^{\frac{3-\alpha-\beta_\gamma}
{\gamma}}t_n^{\sigma-\frac{\beta_\gamma}{\gamma}}+
\tau^{\beta_\gamma}\tau^{2-\alpha-\beta_\gamma}\sum_{j=2}^{n-1}\tau_jt_j^{\sigma-\frac{2-\alpha}
{\gamma}-1}\Big)\nonumber\\
&\leq C_uC_\gamma\Big(\frac{1}{\sigma}+
t_n^{\sigma-\frac{\beta_\gamma}{\gamma}}+
(\rho+1)\gamma t_n^{\sigma-\frac{\beta_\gamma}{\gamma}}\varsigma_{n^*,\gamma}^L\Big)\tau^{\beta_\gamma}\nonumber\\
&\leq C_uC_{\gamma,\rho,T}\Big(\frac{1}{\sigma}+
\varsigma_{n^*,\gamma}^L\Big)\tau^{\beta_\gamma},\label{EQ_0210}
\end{align}
where $\varsigma_{n^*,\gamma}^L$ is defined in \eqref{EQ_etabar} with $n^*:=(\frac{t_n}{\tau_1})^\frac1{\gamma}$ and $C_{\gamma,\rho,T}$ is a constant only depends on $\gamma,\rho,T$.

On the other hand, for $\sigma\in (1,2)$, it follows from \eqref{EQ_1648} that
\begin{align}
\Xi^n_t&\leq C_uC_\gamma\Big(\frac{\tau^{\sigma\gamma}}{\sigma}+\tau^{2-\alpha}\sum_{j=2}^n\tau_j\max_{j\leq k\leq n}t_k^{\frac{\beta_\gamma-2+\alpha}{\gamma}-1}t_k^{\sigma-\frac{\beta_\gamma}{\gamma}}\Big)\nonumber\\
&\leq C_uC_\gamma\Big(\frac{\tau^{\sigma\gamma}}{\sigma}+T^{\sigma-\frac{\beta_\gamma}{\gamma}}\tau^{2-\alpha}
\sum_{j=2}^n\tau_jt_j^{\frac{\beta_\gamma-2+\alpha}{\gamma}-1}\Big).\label{EQ_0330}
\end{align}
For any $\xi\geq 0$, the direct calculation shows that
\begin{align}
\sum_{j=2}^{n}\tau^{j} t_j^{-\xi-1}
\leq \int_{\tau_1}^1s^{-\xi-1}\zd s= \left\{\begin{array}{cc}
\frac{1}{\xi}(\tau_1^{-\xi}-{t_n}^{-\xi}),&\xi>0\\
\ln \frac{t_n}{\tau_1},&\xi= 0.
\end{array}\right.\label{EQ_0335}
\end{align}
Setting $\xi = (2-\alpha-\beta_\gamma)/\gamma$ in \eqref{EQ_0335} and inserting the resulting into  \eqref{EQ_0330}, we have
\begin{align}
\Xi^n_t&\leq C_uC_\gamma\Big(\frac{\tau^{\sigma\gamma}}{\sigma}+\tau^{2-\alpha}\sum_{j=2}^n\tau_j\max_{j\leq k\leq n}t_k^{\frac{\beta_\gamma-2+\alpha}{\gamma}-1}t_k^{\sigma-\frac{\beta_\gamma}{\gamma}}\Big)\nonumber\\
&\leq C_uC_\gamma\Big(\frac{\tau^{\sigma\gamma}}{\sigma}+T^{\sigma-\frac{\beta_\gamma}{\gamma}}\tau^{\beta_\gamma}
\tau_1^{(2-\alpha-\beta_\gamma)/\gamma}\sum_{j=2}^n\tau_jt_j^{\frac{\beta_\gamma-2+\alpha}{\gamma}-1}\Big)\nonumber\\
&\leq C_uC_{\gamma,T}\Big(\frac{1}{\sigma}+
\zeta_{n^*,\gamma}^L\Big)\tau^{\sigma\gamma},\label{EQ_0340}
\end{align}
where $\zeta_{n^*,\gamma}^L$ is defined in \eqref{EQ_1708} with $n^*:=(\frac{t_n}{\tau_1})^\frac1{\gamma}$ and $C_{\gamma,T}$ is a constant only depending on $\gamma,T$.

Thus,  from \eqref{EQ_0210} and \eqref{EQ_0340}, the global consistency error $\Xi^n_t$ on general nonuniform mesh can be bounded by
\begin{equation}\label{EQ_1551}
\begin{aligned}
\Xi^n_t
\leq&  C_uC_{\gamma,\rho,T}\Big(\frac{1}{\sigma}+
\chi_{n^*,\gamma}^L\Big)\tau^{\min\{\sigma\gamma,2-\alpha\}},
\end{aligned}
\end{equation}
where the factor $\chi_{n^*,\gamma}^L$ is defined in \eqref{EQ_cl} with $n^*:=(\frac{t_n}{\tau_1})^\frac1{\gamma}$.

By choosing  $\gamma_{\texttt{opt}}=\max\{1,(2-\alpha)/\sigma\}$, one has the quasi-optimal convergence.

\subsection{Convergence of L1 scheme}
%We now have the convergence of L1 scheme with various meshes.
Inserting \eqref{EQ_alpharobust}, \eqref{EQ_L1uni}, \eqref{gradedconsistency} and \eqref{EQ_1551} into \eqref{EQ_ger3},  we have the following theorem for the $\alpha$-robust convergence of L1 scheme.
\begin{theorem}\label{thm: robustmain}
Suppose the solution $u$ to  \eqref{EQ_equation} satisfies  the regularity assumptions
 \eqref{EQ_2153} and \eqref{SR1}, and
consider the L1 convolution kernels \eqref{L1s} with $\theta=0$.
%numerical scheme \eqref{eq: discrete IBVP1}.
 Let $\kappa_+:=\max\{\kappa,0\}$.
If the maximum step size $\taumax\le1/\sqrt[\alpha]{4\max\{1,\rho\}\Gamma(2-\alpha)\kappa_+}$,
%with the
%maximum step size $\taumax\le1/\sqrt[\alpha]{20\Gamma(2-\alpha)\kappa}$, then
the discrete solution $U_h^n$ is convergent in the $L^2$-norm with $C=(\rho+1)E_\alpha(4\max\{1,\rho\}\kappa_+ t_n^\alpha)$ that
\begin{equation}
\|e_h^n\|\leq C_uC\Big(\|e_h^0\|+
\frac{\tau_1^\sigma}{\sigma}+\sum_{j=2}^n\tau_j\max_{j\leq k\leq n}(t_{k-1}^{\sigma-2}t_k^{\alpha-1}\tau_k^{2-\alpha})+t_n^\alpha h^2\Big), \; 1\leq n\leq N. \label{EQ_0019}
\end{equation}
In particular, if uniform mesh is used, then it holds
\begin{align}%\label{globalError-optimal-simpleScheme}
\mynormb{e_h^n}\leq C_uC
\bra{\|e_h^0\|+(\frac1{\sigma}+\varsigma_{n,1}^L) N^{-\min\{\sigma,2-\alpha\}}
+t_{n}^{\alpha}h^2},  \; 1\leq n\leq N,
\end{align}
if graded mesh is used, then it holds
\begin{align}%\label{globalError-optimal-simpleScheme}
\mynormb{e_h^n}&\leq C_uC
\bra{\|e_h^0\|+(\frac1{\sigma}+\chi_{n,\gamma}^L) N^{-\min\{\sigma\gamma,2-\alpha\}}
+t_{n}^{\alpha}h^2}, \label{EQ_0404}
\end{align}
if \Mss{1} holds, then it holds
\begin{align}%\label{globalError-optimal-simpleScheme}
\mynormb{e_h^n}&\leq C_uC
\bra{\|e_h^0\|+(\frac1{\sigma}+\chi_{n^*,\gamma}^L)\tau^{\min\{\sigma\gamma,2-\alpha\}}
+t_{n}^{\alpha}h^2}, \label{EQ_0405}
\end{align}
where $n^*=(\frac{t_n}{\tau_1})^\frac1{\gamma}$, $\varsigma_{n,\gamma}^L$ and $\chi_{n,\gamma}^L$ are defined by \eqref{EQ_etabar} and \eqref{EQ_cl} respectively.  One can achieve  the optimal convergence by choosing $\gamma_{\text{opt}}=\max\{1,(2-\alpha)/\sigma\}$ in \eqref{EQ_0404}-\eqref{EQ_0405}.
\end{theorem}

\section{Application to Alikhanov's scheme\label{Sec_4}}
We now consider the $\alpha$-robust error estimate of Alikhanov's scheme with the kernels
 $A_0^{(1)}\defeq a_0^{(1)}$, and for $n\geq2$
\begin{equation}\label{eq: weights}
A^{(n)}_{n-k}\defeq\begin{cases}
	a^{(n)}_0+\rho_{n-1}b^{(n)}_1,
	&\text{for $k=n$,}\\
	a^{(n)}_{n-k}+\rho_{k-1}b^{(n)}_{n-k+1}-b^{(n)}_{n-k},
	&\text{for $2\le k\le n-1$,}\\
	a^{(n)}_{n-1}-b^{(n)}_{n-1},
	&\text{for $k=1$,}
\end{cases}
\end{equation}
where the discrete coefficients
$a_{n-k}^{(n)}$ and $b_{n-k}^{(n)}$ $(1\le k\leq n-1)$ are defined as
\begin{align}
&a^{(n)}_0\defeq\frac1{\tau_n}\int_{t_{n-1}}^{t_{n-\theta}}
	\omega_{1-\alpha}(t_{n-\theta}-s)\zd{s}\;\;\mbox{and}\;\; a^{(n)}_{n-k}\defeq\frac{1}{\tau_k}\int_{t_{k-1}}^{t_k}	\omega_{1-\alpha}(t_{n-\theta}-s)\zd{s};\label{eq: an}\\
&b^{(n)}_{n-k}\defeq\frac{2}{\tau_k(\tau_k+\tau_{k+1})}\int_{t_{k-1}}^{t_k}
	(s-t_{k-\frac12})\omega_{1-\alpha}(t_{n-\theta}-s)\zd{s}.  \label{eq: bn}
\end{align}
As proved in \cite{LMZ21}, one can verify the above assumptions hold as follows:
\begin{enumerate}
  \item Assumptions \Ass{1} and \Ass{2} hold with $\pi_A=11/4$,  and  $\rho=7/4$ in \Ass{3};
  \item $\theta = \alpha/2<\theta^{(n)};$
  \item $A^{(n)}_0\leq \frac{24}{11\tau_n}\int_{t_{n-1}}^{t_n} \omega_{1-\alpha}(t_n-s)\leq \frac{24}{11\Gamma(2-\alpha)\tau_n^\alpha}$;
  \item With the regularity assumption \eqref{EQ_2153} and  appropriate regularity in space, the local and history interpolation errors $\Gloc^k$ and $\Ghis^k$ can be bounded by
  \begin{align}
\Gloc^1&\leq C_u\tau_1^\sigma/\sigma\qquad\qquad\qquad \text{and}\quad \Gloc^k\leq C_u t_{k-1}^{\sigma-3}\tau_k^3 ,\label{EQ_2131}\\
\Ghis^1&\leq C_u(\tau_1^\sigma/\sigma+t_1^{\sigma-3}\tau_2^3) \quad \text{and}\quad \Ghis^k\leq C_u (t_{k-1}^{\sigma-3}\tau_k^3+t_k^{\sigma-3}\tau_{k+1}^3).\label{EQ_2132}
\end{align}
\end{enumerate}

Thus, inserting \eqref{EQ_2131} and \eqref{EQ_2132} into \eqref{EQ_Xit}, we have
\begin{align}
  \Xi^n_t:=   C_u\Big(\tau_1^\sigma/\sigma+t_1^{\sigma-3}\tau_2^3+\sum_{j=2}^n \tau_j\max_{j\leq k\leq n}(t_{k-1}^{\sigma-3}t_k^{\alpha-1}\tau_k^{3-\alpha}+t_k^{\sigma+\alpha-4}\tau_{k+1}^3/\tau_k^\alpha)\Big).\label{EQ_2140}
\end{align}
\begin{remark}
Using the similar techniques to Remark \ref{Remark_L1}, a bound of  the global consistency error  presented in \cite{LMZ21} is given by
\begin{equation}
\sum_{k=1}^n P^{(n)}_{n-k} |(R_t)^j|\leq C_u\big(\frac{\tau_1^\sigma}{\sigma}+\frac1{1-\alpha}\max_{2\leq k\leq n}t_k^\alpha t_{k-1}^{\sigma-3}\tau_k^3/\tau_{k-1}^\alpha\big),\label{EQ_1702}
\end{equation}
which also has the factor of $1/(1-\alpha)$, but our new bound \eqref{EQ_2140} does not.
\end{remark}

To investigate the corresponding global consistency error with various meshes, we introduce  two factors which will appear  frequently in the error bounds.  Set
\begin{equation}
\vartheta_\gamma:=\max\{\sigma\gamma,3-\alpha\}.\label{EQ_1552}
\end{equation}
Then we denote by
\begin{equation}
  \varsigma_{n,\gamma}^A:= n^{\alpha-3+\beta_\gamma}\int_{1/n}^{1} s^{\sigma\gamma+\alpha-4} \zd s= \left\{\begin{array}{cc}
\frac{1-n^{\sigma\gamma+\alpha-3}}{3-\alpha-\sigma\gamma},& \sigma\gamma< 3-\alpha,\\
\ln n,&\sigma\gamma=3-\alpha,\\
\frac{1-n^{3-\alpha-\sigma\gamma}}{\sigma\gamma+\alpha-3},&\sigma\gamma>3-\alpha,
\end{array}\right.\label{EQ_va}
\end{equation}
and
\begin{align}
\zeta_{n,\gamma}^A :=n^{\beta_\gamma+\alpha-2}\int_{1/n}^1s^{\beta_\gamma+\alpha-3}\zd s=\left\{\begin{array}{cc}
\frac{1-n^{\sigma\gamma+\alpha-2}}{(2-\alpha)-\sigma\gamma},&\sigma\gamma<2-\alpha,\\
\ln n,&\sigma\gamma\geq 2-\alpha.
\end{array}\right.\label{EQ_za}
\end{align}
For convenience, we further define
\begin{equation}\label{EQ_ca}
  \chi_{n,\gamma}^A:= \left\{\begin{array}{cc}
\varsigma_{n,\gamma}^A,&0<\sigma<1,\\
\zeta_{n,\gamma}^A,&1<\sigma<2.
\end{array}\right.
\end{equation}
From \eqref{EQ_va} and \eqref{EQ_za},  one can see that the factors $\varsigma_{n,\gamma}^A$ and $\zeta_{n,\gamma}^A$ would not blow up for all $\alpha$ and $\sigma$ since
\begin{enumerate}
  \item $\varsigma_{n,\gamma}^A, \zeta_{n,\gamma}^A$  are continuous with respect to $\alpha, \sigma$;
  \item $\varsigma_{n,\gamma}^A\to \ln n$ and $\zeta_{n,\gamma}^A \to \ln n$ as $\sigma+\alpha\to 3$;
      \item $\varsigma_{n,\gamma}^A$ and $\zeta_{n,\gamma}^A$ can be uniformly bounded by $\ln n$ for all $\alpha$ and $\sigma$.
\end{enumerate}
\subsection{The global approximation error $\Xi^n_t$ on uniform mesh}
Taking the constant time step $\tau = T/N$ and $t_k= k\tau$, thus estimate \eqref{EQ_2140} has
\begin{align}
\Xi^n_t&\leq C_u\tau^\sigma(\frac1{\sigma}+1+\sum_{j=2}^nj^{-1}\max_{j\leq k\leq n} ((k-1)^{\sigma-3}k^\alpha+k^{\sigma+\alpha-3}))\nonumber\\
&\leq C_u\tau^\sigma(\frac1{\sigma}+1+\sum_{j=2}^nj^{\sigma+\alpha-4})\nonumber\\
&\leq C_u(\frac{\tau^\sigma}{\sigma}+1+t_n^{\sigma-\vartheta_1}\tau^{\vartheta_1} n^{\vartheta_1+\alpha-3}\int_{1/n}^1s^{\sigma+\alpha-4}\zd s)\nonumber\\
&\leq C_u(\frac{1}{\sigma}+t_n^{\sigma-\vartheta_1} \varsigma_{n,1}^A)\tau^{\vartheta_1}\leq C_uC_T(\frac{1}{\sigma}+ \varsigma_{n,1}^A)\tau^{\vartheta_1},\label{EQ_1614}
%&\leq C_u\tau^\sigma(\frac1{\sigma}+1+n^{\sigma+\alpha-3}\int_{1/n}^{1} s^{\sigma+\alpha-4} \zd s),
\end{align}
where one uses the fact $\frac{k}{k-1}\leq 2$ and $\vartheta_\gamma,  \varsigma_{n,\gamma}^A$ are defined in \eqref{EQ_1552} and \eqref{EQ_va}, respectively.
\subsection{The global approximation error $\Xi^n_t$ on graded mesh}
We now consider the graded mesh $t_k=T\bra{k/N}^{\gamma}$ with $\gamma \geq 1$.
Based on fact $\frac{k}{k-1}\leq 2$ for any $k\geq 2$,  \eqref{EQ_2140} can be reduce to
\begin{align}
  \Xi^n_t&\leq C_uC_{\gamma,T}\Big(\frac1{\sigma} N^{-\sigma\gamma}+\sum_{j=2}^n \tau_j\max_{j\leq k\leq n}( t_k^{\sigma+\alpha-4}\tau_k^{3-\alpha})\Big)\nonumber\\
&\leq C_uC_{\gamma,T}\Big(\frac1{\sigma} N^{-\sigma\gamma}+\sum_{j=2}^n \tau_j\max_{j\leq k\leq n}( k^{(\sigma-1)\gamma-(3-\alpha)}N^{-(\sigma-1)\gamma})\Big). \label{EQ_0055}
\end{align}
Note that $(\sigma-1)\gamma-(3-\alpha)<0$ for any $\gamma\geq 1$ if $\sigma\in (0,1)$. Thus, it holds for $\sigma\in (0,1)$ that
\begin{align}
  \Xi^n_t&\leq C_uC_{\gamma,T}\Big(\frac1{\sigma}N^{-\sigma\gamma} +t_n^{\sigma-\vartheta_\gamma/\gamma } n^{\vartheta_\gamma-(4-\alpha)}N^{-\vartheta_\gamma}+N^{-\sigma\gamma}\sum_{j=2}^{n-1}j^{\sigma\gamma-(4-\alpha)}\Big).
  \label{EQ_0040}
\end{align}
 Setting $\xi = (3-\alpha)/\gamma-\sigma$ in \eqref{EQ_0116} and inserting the resulting into  \eqref{EQ_0040}, we have
 \begin{align}
   \Xi^n_t&\leq C_uC_{\gamma,T}\Big(\frac1{\sigma} +t_n^{\sigma-\vartheta_\gamma/\gamma}\varsigma_{n,\gamma}^A\Big)N^{-\vartheta_\gamma}.\label{EQ_0053}
 \end{align}

 On the other hand, for $\sigma\in(1,2)$, it follows from \eqref{EQ_0055} that
 \begin{align}
  \Xi^n_t
&\leq C_uC_{\gamma,T}\Big(\frac1{\sigma} N^{-\sigma\gamma}+N^{-\vartheta_\gamma}\sum_{j=2}^n j^{-1}\max_{j\leq k\leq n}\big(\frac{k}{N}\big)^{\sigma\gamma-\vartheta_\gamma}k^{\vartheta_\gamma-(3-\alpha)}\Big)\nonumber\\
&\leq C_uC_{\gamma,T}\Big(\frac1{\sigma} +\sum_{j=2}^nj^{\vartheta_\gamma-(4-\alpha)}\Big)N^{-\vartheta_\gamma}\nonumber\\
&\leq C_uC_{\gamma,T}\Big(\frac1{\sigma} +\zeta_{n,\gamma}^A\Big)N^{-\vartheta_\gamma},\label{EQ_1549}
\end{align}
 where the last inequality uses \eqref{EQ_0117} with  $\xi=3-\alpha-\vartheta_\gamma$ and $\zeta_{n,\gamma}^A$ is defined in \eqref{EQ_za}.
Thus,  from \eqref{EQ_0053} and \eqref{EQ_1549}, the global consistency error $\Xi^n_t$ on graded meshes can be bounded as
\begin{equation}\label{EQ_1914}
\begin{aligned}
\Xi^n_t
\leq&  C_uC_{\gamma,T} \braB{\frac{1}{\sigma} +\chi_{n,\gamma}^A }N^{-\min\{\gamma\sigma,3-\alpha\}},
\end{aligned}
\end{equation}
where the factor $\chi_{n,\gamma}^A$  is defined in \eqref{EQ_ca}.
\begin{remark}
In \cite[Lemma 25]{CS21}, the authors consider a case of $\sigma=\alpha$ and present the $\alpha$-robust estimate of the global consistency error as
\begin{align}
\sum_{j=1}^k P^{(k)}_{k-j}\|(R_t)_h^j\|\leq C\frac{11e^{\gamma}\Gamma(1+\ell_N-\alpha)}{4\Gamma(1+\ell_N)}T^\alpha\left(\frac{t_n}{ T}\right)^{\ell_N}N^{-\min\{\alpha \cy{\gamma} ,3-\alpha\}},\label{EQ_cs21}
\end{align}
where $\ell_N :=1/\ln N$.
Note that $\Gamma(1+\ell_N-\alpha)\to \Gamma(\ell_N)$  as $\alpha\to 1^-$, which will tend to infinity as $N\to \infty$ and cause the loss of accuracy. Our estimate  \eqref{EQ_1914} remains the optimal convergence order even if $\alpha\to 1^-$ whenever $\gamma\neq 2$.
\end{remark}
\subsection{The global error $\Xi^n_t$ on general nonuniform mesh}
\indent
We now consider the global error $\Xi^n_t$ on general nonuniform  mesh \Mss{1}.
Under \Mss{1},  \eqref{EQ_2140} can be reduced to
\begin{align}
  \Xi^n_t&\leq    C_uC_\gamma\Big(\tau^{\sigma\gamma}/\sigma+\sum_{j=2}^n \tau_j\max_{j\leq k\leq n}t_k^{\sigma+\alpha-4}\tau_k^{3-\alpha}\Big)\nonumber\\
  &\leq    C_uC_\gamma\Big(\tau^{\sigma\gamma}/\sigma+\tau^{3-\alpha}\sum_{j=2}^n \tau_j\max_{j\leq k\leq n}t_k^{\sigma-1-\frac{3-\alpha}{\gamma}}\Big).\label{EQ_0131}
\end{align}
For $\sigma\in (0,1)$, it is easy to check $ \sigma-1-\frac{3-\alpha}{\gamma}<0$ for all $\gamma\geq 1$. Then it holds for $\sigma\in (0,1)$ that
\begin{align}
  \Xi^n_t
  &\leq    C_uC_\gamma\Big(\tau^{\sigma\gamma}/\sigma+\tau^{3-\alpha}\sum_{j=2}^n \tau_jt_j^{\sigma-1-\frac{3-\alpha}{\gamma}}\Big)\nonumber\\
  &\leq  C_uC_\gamma\Big(\tau^{\sigma\gamma}/\sigma+\tau^{\vartheta_\gamma}\big(\frac{\tau_1}{t_n}\big)^{\frac{3-\alpha-\vartheta_\gamma}
{\gamma}}t_n^{\sigma-\frac{\vartheta_\gamma}{\gamma}}+\tau^{3-\alpha}\sum_{j=2}^{n-1} \tau_jt_j^{\sigma-1-\frac{3-\alpha}{\gamma}}\Big)\nonumber\\
  &\leq  C_uC_{\gamma,\rho,T}\Big(\frac1{\sigma}+t_n^{\sigma-\frac{\vartheta_\gamma}{\gamma}}\varsigma_{n^*,\gamma}^A
  \Big)\tau^{\vartheta_\gamma}\leq C_uC_{\gamma,\rho,T}\Big(\frac1{\sigma}+\varsigma_{n^*,\gamma}^A
  \Big)\tau^{\vartheta_\gamma},\label{EQ_0151}
\end{align}
where ones uses \eqref{EQ_0201} with $\xi = (3-\alpha)/\gamma-\sigma$, $n^*:=(\frac{t_n}{\tau_1})^{\frac1{\gamma}} $ and $\vartheta_\gamma, \varsigma_{n,\gamma}$ are defined in \eqref{EQ_1552} and \eqref{EQ_va}.

On the other hand, for $\sigma\in (1,2)$, it follows from \eqref{EQ_0131} that
\begin{align}
\Xi^n_t&\leq C_uC_\gamma\Big(\frac{\tau^{\sigma\gamma}}{\sigma}+\tau^{3-\alpha}\sum_{j=2}^n\tau_j\max_{j\leq k\leq n}t_k^{\frac{\vartheta_\gamma-3+\alpha}{\gamma}-1}t_k^{\sigma-\frac{\vartheta_\gamma}{\gamma}}\Big)\nonumber\\
&\leq C_uC_\gamma\Big(\frac{\tau^{\sigma\gamma}}{\sigma}+T^{\sigma-\frac{\vartheta_\gamma}{\gamma}}\tau^{3-\alpha}
\sum_{j=2}^n\tau_jt_j^{\frac{\vartheta_\gamma-3+\alpha}{\gamma}-1}\Big)\nonumber\\
&\leq C_uC_{\gamma,T}\Big(\frac{1}{\sigma}+
\zeta_{n^*,\gamma}^A\Big)\tau^{\sigma\gamma},\label{EQ_0207}
\end{align}
where the last inequality uses \eqref{EQ_0335} with $\xi = (3-\alpha-\vartheta_\gamma)/\gamma$, $\zeta_{n^*,\gamma}^A$ is defined in \eqref{EQ_za} with $n^*:=(\frac{t_n}{\tau_1})^\frac1{\gamma}$.

Thus, from \eqref{EQ_0151} and \eqref{EQ_0207},
the global consistency error $\Xi^n_t$ on general nonuniform mesh can be bounded as
\begin{equation}\label{EQ_0213}
\begin{aligned}
\Xi^n_t
\leq&  C_uC_{\gamma,\rho,T}\Big(\frac{1}{\sigma}+
\varsigma_{n^*,\gamma}^A\Big)\tau^{\min\{\sigma\gamma,3-\alpha\}},
\end{aligned}
\end{equation}
where the factor $\chi_{n^*,\gamma}^A$ is defined in \eqref{EQ_ca}  with $n^*:=(\frac{t_n}{\tau_1})^\frac1{\gamma}$.

\subsection{Convergence of Alikhanov's scheme}
We now report the convergence of Alikhanov's scheme.
\begin{theorem}\label{thm: robustmain1}
Suppose that the solution $u$ of \eqref{EQ_equation} has the regularity
assumption \eqref{EQ_2153} and \eqref{SR1}, and
consider the Alikhanov's convolution  kernels \eqref{eq: weights}. Let $\kappa_+:=\max\{\kappa,0\}$ and $\theta = \alpha/2$.
If
%$\theta = \alpha/2$, $\rho=7/4$ and
the maximum step size satisfies $\taumax\le1/\sqrt[\alpha]{11\Gamma(2-\alpha)\kappa+}$,
%with the
%maximum step size $\taumax\le1/\sqrt[\alpha]{20\Gamma(2-\alpha)\kappa}$, then
the discrete solution $U_h^n$ ($1\leq n\leq N$) is convergent in the $L^2$-norm by
\begin{align}
\mynormb{e^n_h}
	\le C_u E_\alpha(20\kappa_+ t_n^\alpha) \Big(&\mynormb{e^0_h}+\frac{\tau_1^\sigma}{\sigma}
	+t_n^\alpha\max_{2\le k\le n}t_{k-1}^{\sigma-2}\tau_k^2 + t_n^\alpha h^2\nonumber\\
	&+ \sum_{j=2}^n \tau_j\max_{j\leq k\leq n}(t_{k-1}^{\sigma-3}t_k^{\alpha-1}\tau_k^{3-\alpha}+t_k^{\sigma+\alpha-4}\tau_{k+1}^3/\tau_k^\alpha)\Big).\label{EQ_0329}
\end{align}
In particular, if uniform mesh is used, then it holds with $C= C_u C_TE_\alpha(20\kappa_+ t_n^\alpha)$
\begin{equation*}
  \mynormb{e^n_h}\leq
 C\bra{\mynormb{e^0_h}+(1/\sigma+\varsigma_{n,1}^A)N^{-\min\{\sigma, 2\}}+t_n^\alpha h^2}, \quad 1\le n\le N,
\end{equation*}
if graded mesh is used, then it holds with $C= C_u C_{\gamma,T} E_\alpha(20\kappa_+ t_n^\alpha)$
\begin{align}
\mynormb{e^n_h}&\leq
C \bra{\mynormb{e^0_h}+(1/\sigma+ \chi_{n,\gamma}^A) N^{-\min\{\sigma\gamma, 2\}}+t_n^\alpha h^2}, \quad 1\le n\le N,
\end{align}
 if \Mss{1} also holds, then it holds with $C= C_u C_{\gamma,\rho,T} E_\alpha(20\kappa_+ t_n^\alpha)$
\begin{align}
\mynormb{e^n_h}&\leq
C \bra{\mynormb{e^0_h}+(1/\sigma+\chi_{n^*,\gamma}^A) N^{-\min\{\sigma\gamma, 2\}}+t_n^\alpha h^2}, \quad 1\le n\le N,
\end{align}
where $n^*=(\frac{t_n}{\tau_1})^\frac1{\gamma}$, $\varsigma_{n,\gamma}^A$ and $\chi_{n,\gamma}^A$ are defined by \eqref{EQ_va} and \eqref{EQ_ca} respectively.
\end{theorem}

It is remarkable that if the assumption~$\Mss{1}$ holds, we have
$\tau_1\le C_{\gamma}\tau^{\gamma}$ and
$$t_{k}^{\alpha}t_{k-1}^{\sigma-3}\tau_k^{3}/\tau_{k-1}^{\alpha}
\le C_{\gamma}t_k^{\sigma-(3-\alpha)/\gamma}\tau^{3-\alpha},$$
which implies that the Alikhanov formula $\dfd{\alpha}v^{n-\theta}$ approximates
${}^{\text{C}}_0\fd{\alpha}u(t_{n-\theta})$ in the order of $O(\tau^{3-\alpha})$ if
$\gamma\ge(3-\alpha)/\sigma$. But, noting the difference~$u(t_{n-\theta})-u^{n-\theta}$ in time,
the convergence order  is still limited to the order of~$O(\tau^{2})$.
%\begin{remark}
%In \cite{LMZ21}, the error estimate is given with $C=C_uE_\alpha(20\kappa_+ t_n^\alpha)$ as
%\begin{align}
%\|e_h^n\|\leq C\Big(\frac{\tau_1^\sigma}{\sigma}+\frac1{1-\alpha}\max_{2\leq k\leq n}t_k^\alpha t_{k-1}^{\sigma-3}\tau_k^3/\tau_{k-1}^\alpha+t_n^\alpha\max_{2\leq k\leq n}t_{k-1}^{\sigma-2}\tau_k^2+t_n^\alpha h^2
%\Big).\label{EQ_esAl}
%\end{align}
%The error bound is  $\alpha$-nonrobust since it contains a factor $1/(1-\alpha)$ and will blow up as $\alpha\to 1^-$. The factor is circumvented in  our $\alpha$-robust error bound  \eqref{EQ_0329} by  a summation-type bound $ \sum_{j=2}^n \tau_jt_j^{-1}\max_{j\leq k\leq n}(t_{k-1}^{\sigma-3}t_k^\alpha\tau_k^{3-\alpha}+t_k^{\sigma+\alpha-3}\tau_{k+1}^3/\tau_k^\alpha)$ in \eqref{EQ_0329}.
%Actually, \eqref{EQ_esAl} can be reproduced without the factor $1/(1-\alpha)$ by \eqref{EQ_0329} with  a further estimate of the  summation-type bound
%\begin{align*}
% \sum_{j=2}^n \tau_jt_j^{-1}&\max_{j\leq k\leq n}(t_{k-1}^{\sigma-3}t_k^\alpha\tau_k^{3-\alpha}+t_k^{\sigma+\alpha-3}\tau_{k+1}^3/\tau_k^\alpha)\\
% &\leq \ln \frac{t_n}{\tau_1} \max_{2\leq k\leq n}(t_{k-1}^{\sigma-3}t_k^\alpha\tau_k^{3-\alpha}+t_k^{\sigma+\alpha-3}\tau_{k+1}^3/\tau_k^\alpha)
%\end{align*}
%\end{remark}
\section{Conclusion\label{Sec_5}}
A novel $\alpha$-robust error analysis for convolution-type schemes with general nonuniform time-step is developed for reaction-subdiffusion equations. The new analysis not only circumvents the factor $1/(1-\alpha)$ and hence avoids the factor blow-up phenomenon, but also produces the optimal convergence order for any given $\alpha$ and $\sigma$. While  $\alpha\to 1^{-}$, a factor $\ln n$ appears in the error bound, which is acceptable for quasi-optimal convergence order.
We would like to emphasize that our analysis framework can be extended naturally to nonlinear fractional problems.
%As a final remark, the analysis in this paper is aimed at the L2 error estimate, the framework is also suitable for the H1 error estimate once the scheme satisfies assumptions {\bf A1-A3}.

\bibliographystyle{plain}
\bibliography{AR}
\end{document}